\documentclass[12pt,a4paper]{article}

\usepackage{amsfonts}
\usepackage{amssymb}
\usepackage{amsmath}
\usepackage{amsthm}
\usepackage{url}
\usepackage{enumerate}
\usepackage{hyperref}
\def\IC{\mathbb{C}}
\def\IR{\mathbb{R}}

\def\El{\mathcal{E}\ell}
\def\BH{\mathcal{B}(H)}
\def\trace{\mathop{\mbox{trace}}}
\def\tgm{\mathop{\mbox{tgm}}}
\def\hatted{{\widehat{\phantom{m}}}}
\numberwithin{equation}{section}

\newtheorem{proposition}{Proposition}[section]
\newtheorem{lemma}[proposition]{Lemma}
\newtheorem{theorem}[proposition]{Theorem}
\newtheorem{corollary}[proposition]{Corollary}
\newtheorem{remarks}[proposition]{Remarks}

\newtheorem{definition}[proposition]{Definition}
\newtheorem{example}[proposition]{Example}
\newtheorem{notation}[proposition]{Notation}

\begin{document}
\title{Some Formulae for Norms of Elementary Operators}
\author{Richard M. Timoney}
\maketitle

\begin{abstract}
We present a formula for the norm of an elementary operator on a
$C^*$-algebra that seems to be new. The formula involves (matrix)
numerical ranges and a kind of geometrical mean for positive matrices,
the tracial geometric mean, which seems not to have been studied
previously and has interesting properties. In addition, we characterise
compactness of elementary operators.
\end{abstract}

We consider an elementary operator $Tx = \sum_{j=1}^\ell a_j x b_j$
(with $x \in A$ a $C^*$-algebra and $a_j, b_j \in M(A)$, with $M(A)$
the multiplier
algebra of $A$).  We denote the class of elementary operators $T
\colon A \to A$ by
$\El(A)$. Specifically, we address the question of finding
a concrete formula for the operator norm $\|T\|$. This problem has been
considered (at least implicitly) over a long period by several
authors and there are solutions known under various special
circumstances (generalised derivations \cite{Stampfli}, antiliminal
by abelian $C^*$-algebras \cite{AMS}, for example). See \cite{MMNormProb}
for a recent survey of the problem, or see \cite[\S5.4]{AraMatBk}
for a brief summary of its importance. One way to view the literature
that relates to the problem is to separate two strands
of problems. One strand concentrates on elementary operators of a rather
special form (with $\ell \leq 2$) and the other (for arbitrary
$\ell$) has relied largely on dealing with the completely bounded norm
$\|T\|_{cb}$ and the Haagerup tensor norm estimate
$\|T\|_{cb} \leq \left\| \sum_{j=1}^\ell a_j \otimes b_j \right\|_h$.

For special forms where $\ell\leq 2$, the case $\ell=1$ is well
understood (see \cite{MMNormProb}). There is a
significant body of literature dealing with (inner) derivations
$\delta_a(x) = ax - xa$ and the estimate $\|\delta_a\| \leq 2 \inf
\|a - z\|$ with the infimum over $z$ in the centre $Z(M(A))$ of
$M(A)$ (see references in \cite[\S4.1, \S4.6]{AraMatBk} and
\cite{MMNormProb}).  In the case $A= \BH$ is the algebra of all
bounded linear operators on a Hilbert space $H$ (or $A=\mathcal{K}(H)$,
the compacts) $Z(M(A))$ is just scalar multiples of the identity and
\cite{Stampfli} showed equality in this estimate for $\|\delta_a\|$.
Subsequent work has generalised this equality to various classes
of $C^*$-algebras but \cite[3.2, 3.3]{Somerset94} implies a
characterisation of those $A$ where equality always holds (those
where all Glimm ideals of $M(A)$ are 3-primal). Moreover in case
this condition is not true, then there is $a \in M(A)$ with
$\|\delta_a\| \leq \sqrt{3} \inf_{z \in Z(M(A))} \|a - z\|$ (and
further related work is to be found in
\cite{Arch78,Somerset93,Somerset94,Somerset97,ArchSomerset04}).
An example of \cite{AST}
shows that the condition on Glimm ideals of $M(A)$ is difficult to
relate to the structure of the primitive ideal space of $A$, so
that the results are perhaps most satisfactory in the unital case
where $M(A) = A$. An alternative approach is given in
\cite[4.1.23]{AraMatBk} where it is shown that (for general $A$)
$\|\delta_a\| = 2 \inf\{\|a - z\| :z \in Z({}^cM(A)) \}$ with
${}^cM(A)$ the bounded central closure of $M(A)$.

For generalised (inner) derivations $\delta_{a,b}(x) = ax - xb$,
there are results that are rather less comprehensive than for
$\delta_a$. In particular \cite{Stampfli} shows $\|\delta_{a,b}\|
= \inf_{\lambda \in \IC} \|a - \lambda\| + \|b - \lambda\|$ when
$A = \BH$ and there is a result in terms of representations of
${}^cM(A)$ and $Z({}^cM(A))$ in  \cite[4.1.23]{AraMatBk}. One may invoke
operator space methods and the fact that $\|\delta_{a,b}\| =
\|\delta_{a,b}\|_{cb}$ as another approach. The `obvious' estimate
that arises from taking account of the centre is then $\|\delta_{a,b}\|
\leq \|a \otimes 1 - 1 \otimes b\|_{Z,h}$ where $\| \cdot \|_{Z,h}$
is the central Haagerup tensor norm on $M(A) \otimes M(A)$. Results
of \cite{AST} show that there is equality in this estimate (for all
$a, b \in M(A)$) when all Glimm ideals of $M(A)$ are 5-primal, but
that it is not sufficient for all Glimm ideals to be 3-primal.

The case of Jordan mappings $J_{a,b}(x) = axb + bxa$ was also the
subject of several papers and a conjecture of M. Mathieu that
$\|J_{a,b}\| \geq \|a\| \|b\|$ when $A = \BH$ (or more generally when
$A$ is a prime $C^*$-algebra) has recently been proved to be true. In
\cite{MagTurn} a weaker result $\|J_{a,b}\|_{cb} \geq \|a\| \|b\|$ was
shown and in \cite{Timoneyaxb} a different proof of this was given
along with a proof of the conjecture using results from
\cite{TimoneyNE}. The conjecture was also shown slightly earlier
in \cite{BBR} by quite different methods.

Turning to progress on general elementary operators (where $\ell$ can
be large) the most satisfactory progress to date deals with
$\|T\|_{cb}$ rather than $\|T\|$. It is shown in \cite{AMS} that
$\|T\|_{cb} = \|T\|$ for all elementary operators $T$ if and only if
the $C^*$ algebra $A$ is `antiliminal by abelian'. It is shown in
\cite{Somerset98} that $\|T\|_{cb} = \left\|\sum_{j=1}^\ell a_j \otimes b_j
\right\|_{Z,h}$ for each elementary operator
$Tx = \sum_{j=1}^\ell a_j x b_j$ on $A$
if all
Glimm ideals of $M(A)$ are primal, and \cite{AST} establishes the
converse. In fact \cite{AST} also gives a necessary and sufficient
condition for this equality when restricted to specific $\ell$ (all
Glimm ideals of $M(A)$ should be $N$-primal for $N = \ell^2 + 1$).

Note that a result of \cite{Magajna93} says that elementary operators
can approximate arbitrary bounded linear operators on
$A$ which preserve all ideals of $A$ (albeit in the strong operator
topology). It follows that in many ways elementary operators must be
typical, especially if $A$ is simple. For general $A$, the special
properties of elementary operators can be captured via the concepts of
central bimodule homomorphisms used in \cite[\S5.3]{AraMatBk}.

In this work,
we give slightly different approaches to our formula for the norm
$\|T\|$ of an elementary operator $T$ in Theorems
\ref{thmformulaA} and \ref{ThmFormulaS1} for the case $A = \BH$ and
in Theorem~\ref{CStarFormula} for the case of general $A$.  The
formulae involve the tracial geometric mean as defined in
Definition~\ref{DEFtgm} and we give some basic properties of this
mean in Remarks~\ref{tgmremarks}. The formulae also involve matrix-valued
numerical ranges considered in \cite{TimoneyNE} but are independent
of the different possible expressions for the elementary operator.
We use the formulae to establish some bounds on the growth of
$k$-norms of $T \in \El(A)$ as $k$ increases and exploit continuity in
our formula to give a new proof
of a result of \cite{AMS} showing that there is no growth at all
if $A$ is antiliminal (Corollary~\ref{antiliminalcoroll}).

The question of compactness of elementary operators has also been
investigated under several circumstances (for example on the Calkin
algebra --- see references given in \cite[p. 232]{AraMatBk}) but a
characterisation for general $C^*$-algebras has eluded proof.
In Theorem~\ref{CompactnessThm} we characterise compactness of $T
\in \El(A)$ via the possibility of choosing compact $a_j$ and $b_j$.

This work arose out of discussions with Aleksej Turn\v sek in Dublin
in January 2004 and we thank him for several helpful comments on
earlier drafts of this paper.
Thanks also to the referee for several suggestions.

\section{The case $\BH$}
For the moment we take $A = \mathcal{B}(H) = M(A)$ 
and establish a formula for $\|T\|$ in
that case. Later we will extend the formula to general $A$.

Recall again
the upper bound $\|T\| \leq \left\| \sum_{j=1}^\ell
a_j \otimes b_j \right\|_h$ in terms of the Haagerup norm $\| \cdot
\|_h$ on $\BH \otimes \BH$ (see \cite[5.4.7]{AraMatBk} for example).
We know that equality holds in case the operators $a_i a_j^*$ commute
and the operators $b_j^* b_k$ also commute \cite[Theorem 3.3 and
Remark 2.5]{TimoneyNE}. In general the inequality is strict.

For $\eta, \xi \in H$
we use the notation $\eta \otimes \xi^*$
for the rank one operator on $H$ with
$(\eta \otimes \xi^*)(\theta) = \langle \theta, \xi \rangle \eta$.

\begin{lemma}
\label{firstlemma}
For $T \in \El(\BH)$, $Tx = \sum_{j=1}^\ell a_j x b_j$,
we have
\[
\|T\| =
\sup_{p_1, p_2} 
\left\| \sum_{j=1}^\ell (p_1 a_j) \otimes (b_j
p_2) \right\|_h
\]
where $p_1, p_2 \in \BH$ are rank one projections
($p_i^2 = p_i = p_i^*$ ($i = 1, 2$)).
\end{lemma}

\begin{proof}
Let $p_1 = \xi \otimes \xi^*$ and $p_2 = \eta \otimes \eta^*$
be
1-dimensional projections (where $\eta, \xi \in H$ are unit vectors).
We look at the operator 
\[
T_{p_1, p_2}(x) = \sum_{j=1}^\ell p_1 a_j x b_j p_2,
\]
an operator with a 1-dimensional range. Specifically it is the operator 
\[
x
\mapsto \langle (Tx) \eta, \xi \rangle \xi \otimes \eta^*
\]
and thus almost a linear functional.

For this operator, $(p_1a_i)(p_1a_j)^*$ are commuting and so are
$(b_ip_2)^*(b_jp_2)$. Hence
\[
\| T_{p_1, p_2} \| = \left\| \sum_{j=1}^\ell (p_1a_j) \otimes (b_j
p_2) \right\|_h
\]
by the remarks above (\cite[Theorem 3.3 and Remark
2.5]{TimoneyNE}). Alternatively, one can appeal to the fact that the
norm of a linear functional is the same as its completely bounded
norm, hence $\| T_{p_1, p_2} \| = \| T_{p_1, p_2} \|_{cb} =$ the
Haagerup tensor norm (by a result of Haagerup --- see
\cite[5.4.7]{AraMatBk} for example).

Now, clearly
\begin{eqnarray*}
\|T\| &=& \sup \{ \|Tx\| : x \in \mathcal{B}(H), \|x\| \leq 1\}
\\
&=& \sup \{ \Re \langle (Tx)\eta, \xi \rangle :\\
&& \qquad x \in \mathcal{B}(H), \|x\| \leq 1,
\xi, \eta \in H, \|\xi\| = \|\eta\| = 1\}
\\
&=& \sup \{ \Re \langle (T_{p_1, p_2}x)\eta, \xi
\rangle : x , \xi, \eta \mbox{ as above} \}\\
&\leq& \sup_{p_1, p_2} \| T_{p_1, p_2}\|
\end{eqnarray*}

Since $\| T_{p_1, p_2}\| \leq \|T\|$, the lemma follows.
\end{proof}

\begin{notation}
\rm
Following the ideas in \cite{TimoneyNE}, we introduce some notation
relating to matrix numerical ranges.

For $\mathbf{b} = [b_1, b_2, \ldots, b_\ell]^t$ a column of operators
$b_j \in \BH$ we consider the matrix of operators
\[
Q(\mathbf{b}) = ( b_i^* b_j)_{i, j = 1}^\ell
= [b_1, b_2, \ldots, b_\ell]^*[b_1, b_2, \ldots, b_\ell]
\in M_\ell(\BH) \equiv \mathcal{B}(H^\ell)
\]
and for $\eta \in H$
\[
Q(\mathbf{b}, \eta) = ( \langle b_i^* b_j \eta, \eta \rangle )_{i, j =
1}^\ell
= ( \langle b_j \eta, b_i \eta \rangle )_{i, j = 1}^\ell
\in M_\ell(\IC).
\]

A matrix numerical range associated with $\mathbf{b}$ that was
considered in \cite{TimoneyNE} is
\[
W_m(\mathbf{b}) = \{ Q(\mathbf{b}, \eta)^t : \eta \in H, \|\eta\| = 1 \}
\subset M_\ell(\IC).
\]
$W_{m,e}(\mathbf{b})$
denotes the set of matrices
of maximal trace
in
the closure of $W_m(\mathbf{b})$.
Recall that the maximal trace is in fact
$\|\mathbf{b}\|^2 = \left\| \sum_{j=1}^\ell b_j^* b_j \right\|$.

In \cite{TimoneyNE} it is shown that there is equality in the Haagerup
estimate for the norm of the elementary operator $T$,
\[
\|T\| \leq \| \mathbf{a} \| \|\mathbf{b}\| = \sqrt{
\left\| \sum_{j=1}^\ell a_j a_j^* \right\|
\left\| \sum_{j=1}^\ell b_j^* b_j \right\|}
\leq
\frac12 \left(
\| \mathbf{a} \|^2 + \|\mathbf{b}\|^2
\right)
\]
if and only if $W_{m,e}(\mathbf{a}^*) \cap W_{m,e}(\mathbf{b}) \neq
\emptyset$ (for $\mathbf{a} = [a_1, a_2, \ldots, a_\ell] \in \BH^\ell$
a row
matrix of operators $a_j \in \BH$ and $\mathbf{a}^* = [a_1^*, a_2^*,
\ldots, a_\ell^*]^t$ a column).
\end{notation}

\begin{definition}
\label{DEFtgm}
For
two positive semidefinite $\ell
\times \ell$ matrices $X$ and $Y$
we define the
\emph{tracial geometric mean} of
$X$ and $Y$
by
\[
\tgm(X, Y) = \trace \sqrt{\sqrt{X}Y\sqrt{X}}
\]
(where, of course, the square roots
mean the positive semidefinite square roots).
\end{definition}

Here is one version of our formula for the norm of an elementary
operator.

\begin{theorem}
\label{thmformulaA}
For $\mathbf{a} = [a_1, a_2, \ldots, a_\ell] \in \BH^\ell$ (a row
matrix of operators $a_j \in \BH$) and
$\mathbf{b} = [b_1, b_2, \ldots, b_\ell]^t \in \BH^\ell$
a column, and $Tx = \sum_{j=1}^\ell a_j x b_j$ an elementary operator,
we have
\[
\|T\| = \sup \left\{
\tgm \left( Q(\mathbf{a}^*, \xi), Q(\mathbf{b}, \eta)
\right)
: \xi, \eta \in H, \|\xi\| = \|\eta\| = 1
\right\}
\]
\end{theorem}

\begin{proof}
Note first that the result will follow from the lemma once we establish that
\[
\|T_{p_1, p_2}\| = 
\left\| \sum_{j=1}^\ell (p_1 a_j) \otimes (b_j
p_2) \right\|_h =
\trace \sqrt{
Q(\mathbf{a}^*, \xi)
}
Q(\mathbf{b}, \eta)
\sqrt{
Q(\mathbf{a}^*, \xi)
}
\]
when $p_1 = \xi \otimes \xi^*$, $p_2 = \eta \otimes \eta^*$.

We now fix unit vectors $\xi, \eta \in H$. Notice that the matrix
numerical range of the tuple $\mathbf{b} {p_2}=
[b_1 p_2, b_2 p_2, \ldots , b_\ell
p_2]^t$ consists of 
\[
\{ |\lambda|^2 Q(\mathbf{b} {p_2}, \eta)^t :
\lambda \in \IC, |\lambda| \leq 1 \}
\]
(if $H$ is not one-dimensional; for the trivial
one-dimensional case we only have $|\lambda|=1$). Moreover
$Q(\mathbf{b} {p_2}, \eta) = Q(\mathbf{b}, \eta)$ for $\mathbf{b} =
[b_1, b_2, \ldots, b_\ell]^t$. Similar remarks
apply to the tuple $\mathbf{a}^*$ and $(p_1 \mathbf{a})^* = [p_1 a_1,
p_1a_2, \ldots, p_1 a_\ell]^*$.
Thus the condition (from \cite{TimoneyNE})
for equality in the Haagerup estimate
\[
\|T_{p_1, p_2}\| \leq \frac{1}{2}\left( \|p_1\mathbf{a}\|^2 +
\|\mathbf{b} {p_2}\|^2 \right)
\]
is
$
Q((p_1\mathbf{a})^*, \xi)^t = Q(\mathbf{b} {p_2}, \eta)^t
$.

We now show how to rewrite $\sum_{j=1}^\ell (p_1a_j) \otimes (b_jp_2)$
so as to get this equality condition satisfied. To simplify the
notation, we assume now that $p_1a_j = a_j$ and $b_j p_2 = b_j$ for $1
\leq j \leq \ell$. Introduce the notation
\[
\mathbf{a} \odot \mathbf{b} = \sum_{j=1}^\ell a_j \otimes b_j
\]
for rows $\mathbf{a}$ and columns $\mathbf{b}$.

Suppose first that the tuples $(a_j)_{j=1}^\ell$ and
$(b_j)_{j=1}^\ell$ are each linearly independent.
Then the different
ways to rewrite $\mathbf{a} \odot \mathbf{b}$ as a sum of the same
type, without going outside the spans of the $a_j$ and $b_j$ (or
equivalently sticking to linearly independent tuples) are all
of the form
\[
(\mathbf{a}\alpha) \odot (\alpha^{-1} \mathbf{b})
\]
for invertible scalar matrices $\alpha \in M_\ell(\IC)$. We have
\[
Q( (\mathbf{a}\alpha)^*, \xi )^t = \alpha^* Q(\mathbf{a}^*, \xi)^t
\alpha = \alpha^* (\langle a_i^* \eta, a_j^* \eta \rangle )_{i, j=1}^\ell \alpha
\]
and
\[
Q(\alpha^{-1} \mathbf{b}, \eta)^t = \alpha^{-1} Q(\mathbf{b}, \eta)^t
\alpha
= \alpha^{-1} (\langle b_i \eta, b_j \eta \rangle )_{i, j=1}^\ell
(\alpha^{-1})^*
\]
For
$Q(\mathbf{a}^*, \xi)$
and
$Q(\mathbf{b}, \eta)$
invertible,
we take 
$\alpha = \alpha_0 \alpha_1$,
\[
\alpha_0 =
\left(\sqrt{Q(\mathbf{a}^*, \xi)^t}\right)^{-1},
\alpha_1 = \left(\sqrt{Q(\mathbf{a}^*, \xi)^t}
Q(\mathbf{b}, \eta)^t \sqrt{Q(\mathbf{a}^*, \xi)^t}\right)^{1/4}.
\]
The effect of $\alpha_0$ is to make $Q( (\mathbf{a}\alpha_0)^*, \xi )$
the identity matrix, and then $\alpha$ is designed so that
we get an equality
\begin{equation}
\label{formul}
Q( (\mathbf{a}\alpha)^*, \xi )^t =
Q(\alpha^{-1} \mathbf{b}, \eta)^t = 
\left(\sqrt{\sqrt{Q(\mathbf{a}^*, \xi)}
Q(\mathbf{b}, \eta) \sqrt{Q(\mathbf{a}^*, \xi)}}\right)^t
.
\end{equation}
Taking $x$ to be a unitary operator in $\BH$ with $x((\alpha^{-1}
\mathbf{b})_j \eta) = (\mathbf{a}\alpha)^*_j \xi$ for $1 \leq j \leq
\ell$, we get $
\langle T_{p_1,p_2}(x)\eta, \xi \rangle =
\langle T(x)\eta, \xi \rangle =$ the trace of the
matrix (\ref{formul}).
Hence we have shown equality in the Haagerup estimate
for $\|T_{p_1, p_2}\|$ arising from $(\mathbf{a}\alpha) \odot
(\alpha^{-1} \mathbf{b})$. Thus $\|T_{p_1, p_2}\|$ is the trace of the
matrix (\ref{formul}).

We have made a linear independence assumption and the calculation
we made requires $Q(\mathbf{a}^*, \xi)$ and $Q(\mathbf{b}, \eta)$ to
be invertible. Though we could perhaps manage without these assumptions and
make use of various notions of generalised inverses, it is
easier to use an approximation argument to deduce the general case. We
embed $\BH$ in $\mathcal{B}(H \oplus H_0)$ where $H_0$ is an
auxiliary Hilbert space of dimension $\ell$. Then we can modify $a_j$
so that $a_j^*\xi$ acquire small mutually
orthogonal contributions in
$H_0$. Similarly for $b_j$ and $b_j \eta$.
This will ensure the assumptions are
valid and then we can take limits as the modifications of $a_j$ and
$b_j$ tend to 0.
\end{proof}

\begin{example}
\rm
If $a_j^*$ have orthogonal ranges and $b_j$ have orthogonal ranges,
then $Q(\mathbf{a}^*, \xi)$ and $Q(\mathbf{b}, \eta)$ are
diagonal matrices and the formula from Theorem~\ref{thmformulaA}
becomes
\[
\|T\| = \sup \sum_{j=1}^\ell \|a_j^*\xi\| \|b_j \eta\|
\]

\end{example}

\begin{remarks}
\label{tgmremarks}
\rm
\begin{enumerate}[(i)]
\item
\label{generalXY}
In Theorem~\ref{thmformulaA}, we are considering the
tracial geometric mean of two positive semidefinite matrices of the
form
$X = Q(\mathbf{a}^*, \xi)$ and $Y = Q(\mathbf{b}, \eta)$.
Such matrices $X$ and $Y$ can be arbitrary
elements of $M_\ell^+(\IC)$ (as long as the dimension of $H$ is not
smaller than $\ell$).

For example, starting with $Y \in
M_\ell^+(\IC)$ we could take $\eta_1, \eta_2, \ldots, \eta_\ell \in
\IC^\ell$ to be the rows of $\sqrt{Y^t}$. Then take $b_j = \eta_j
\otimes \eta^*$ (with $\eta$ any unit vector in the $\ell$-dimensional
Hilbert space $\IC^\ell$) to get $b_j \eta = \eta_j$ and
\[
( \langle b_i \eta, b_j \eta \rangle )_{i, j=1}^\ell
= Y^t \Rightarrow Q(\mathbf{b}, \eta) = Y.
\]
Similarly we can find $a_j = \eta \otimes  \xi_j^*$ so that $X =
Q(\mathbf{a}^*, \eta) = ( \langle a_j^* \eta, a_i^* \eta \rangle )_{i,
j=1}^\ell$ and for $p=p_1 = p_2=\eta \otimes \eta^*$ the operator $T \in
\El(M_\ell)$ with $Tx = \sum_{j=1}^\ell a_j x b_j$ has $T = T_{p,p}$,
$\|T\| = \tgm(X,Y)$.

\item
\label{eigenvaluesversion}
A useful fact
fact is that the eigenvalues of $\sqrt{X}Y\sqrt{X}$ are the same as
those of $\sqrt{X}(\sqrt{X}Y) = XY$. (This well-known fact follows by
conjugation if $\sqrt{X}$ is invertible.)
So if $\lambda_i(XY)$ ($1 \leq i \leq \ell$) denote the
eigenvalues of $XY$ (arranged in non-increasing order, say), then
\[
\tgm(X, Y) = \sum_{i=1}^\ell \sqrt{\lambda_i(XY)}.
\]

\item The symmetry property $\tgm(X, Y) = \tgm(Y, X)$
follows because\\$\lambda_i(XY) = \lambda_i(YX)$. (It can also be
established using $\|T\| = \|T^*\|$, where $T^*$ means $T^*(x) = T(x^*)^*$,
and (\ref{generalXY}).)

\item $\tgm(X, Y)$ has some other desirable properties for
geometric means. If $X = Y$, then $\tgm(X, Y) = \trace X$. If
$X = \lambda I_\ell$ is a multiple of the identity matrix, then
$\tgm(X, Y) = \sqrt{\lambda} \, \trace \sqrt{Y}$.

From the Haagerup estimate for $\|T_{p_1, p_2}\|$ and remark
(\ref{generalXY}) above, we can see that
\[
\tgm(X, Y) \leq \sqrt{(\trace X)(\trace Y)} \leq \frac{\trace X
+ \trace Y}{2}
\]
holds. This is a tracial version of an arithmetic-geometric mean
inequality.

According to a criterion established in \cite{TimoneyNE}
for overall equality in this estimate, equality will only hold if $X =
Y$. It follows that
\[
\tgm(X, Y) = \sqrt{(\trace X)(\trace Y)}
\]
holds only when $X$ and $Y$ are linearly dependent.

\item It is clear that if $u \in M_\ell(\IC)$ is unitary, then
$\tgm(u^*Xu,u^*Yu) = \tgm(X,Y)$ and it follows from
the earlier remark (\ref{eigenvaluesversion}) (or the proof of
Theorem~\ref{thmformulaA})
that if $\alpha \in M_\ell$ is invertible then
\[
\tgm(\alpha^*X\alpha,(\alpha^{-1})Y(\alpha^{-1})^*)
= \tgm(X,Y)
\]

\item It is easy to see that $\tgm(X,Y)$ is monotone in $Y$ and
by symmetry also in $X$. So if $X \leq X_1$ and $Y \leq Y_1$ then
$\tgm(X,Y) \leq \tgm(X_1, Y_1)$.

\item There are various notions of geometric mean for two positive
semidefinite matrices in
the literature (see for example \cite{AndoEtAl}).
One of these is usually denoted $X \# Y$
and can be defined for the positive definite case by
\[
X \# Y = \sqrt{X} \left( X^{-1/2} Y X^{-1/2} \right)^{1/2} \sqrt{X}.
\]
In general $\trace (X \# Y ) \leq
\tgm(X, Y)$ and strict inequality is possible.

T. Ando has kindly provided the following proof of the
inequality (private correspondence). Assuming $X$ and $Y$ are
non-singular, by \cite[D4]{AndoEtAl} $X \# Y = X^{1/2}UY^{1/2}$ for
some unitary $U$. Hence $\trace (X \# Y ) = \trace (UY^{1/2} X^{1/2})$
is at most the trace class
norm
\[
\| Y^{1/2} X^{1/2}\|_1 = 
\trace  (X^{1/2}Y^{1/2}Y^{1/2} X^{1/2}) ^{1/2} = \tgm(X,Y).
\]
The case of general positive
semidefinite $X$ and $Y$ follows by continuity.

To illustrate strict inequality, let
\[
X = \left( \begin{array}{cc}
1 &\frac12\\
\frac12 & 1 \end{array} \right)^2 = 
\left( \begin{array}{cc}
\frac54 & 1\\
1 & \frac54  \end{array} \right), \quad Y = \left( \begin{array}{cc}
1 & 0\\
0 & 0 \end{array} \right).
\]
Then
\[
X \# Y =\left( \begin{array}{cc} \frac{3}{2 \sqrt{5}} & 0\\
0 & 0 \end{array} \right), 
\quad
\tgm(X, Y) = \frac{\sqrt{5}}{2}
\]

\item If $P$ is an orthogonal projection in $M_\ell$ (that is $P = P^2
= P^*$) and $Q = I_\ell - P$,
we define a pinching map by $\mathcal{C}(X) = PXP + QXQ$ and then we
have
\begin{eqnarray*}
\tgm(X, Y)
&\leq& \tgm(\mathcal{C}(X), \mathcal{C}(Y))\\
&=& \tgm(PXP, PYP) + \tgm(QXQ, QYQ)
\end{eqnarray*}

To verify this, we can assume (by a unitary change of basis) that $P$
is the projection of $\IC^\ell$ onto the first $k$ coordinates ($0
\leq k \leq \ell$), in other words a diagonal matrix with $k$ 1's and
$\ell -k$ zeros.
Choose $T \in \El(M_\ell)$ as indicated in remark (\ref{generalXY}) so
that $T = T_{p,p}$ and $\|T\| = \tgm(X,Y)$.
Taking $T_1x = \sum_{j=1}^k a_j x b_j$ and
$T_2 x = \sum_{j=k+1}^\ell a_j x b_j$, we find
\[
\|T\| \leq \|T_1\| + \|T_2\| = \tgm(PXP, PYP) + \tgm(QXQ, QYQ)
\]
when we calculate $\|T_1\|$ and $\|T_2\|$ in a similar way.

\item The maps $X \mapsto \tgm(X,Y) = \tgm(Y,X)$
(with $Y \in M_\ell^+$ fixed) 
satisfy subadditivity
\[
\tgm\left (\sum_{j=1}^n X_j , Y \right) \leq \sum_{j=1}^n \tgm(X_j, Y)
\]
(This can be shown by
using $\tgm(X, Y) = \tgm(\sqrt{Y}X\sqrt{Y}, I_\ell)$
to reduce to the case $Y=I_\ell$
and using a known property of the Schatten $1/2$ quasinorm --- see
\cite[Rotfel'd Theorem IV.2.14]{Bhatia}).
Hence, for $X_j , Y_k \in M_\ell^+$,
\begin{equation}
\label{quasiconvexity}
\tgm \left (\sum_{j=1}^n X_j , \sum_{k=1}^m Y_k \right) \leq \sum_{1
\leq j \leq n, 1 \leq k \leq m} \tgm(X_i, Y_j).
\end{equation}
It follows by Cauchy-Schwarz that
\begin{equation}
\label{quasiconvexityalso}
\tgm \left (\sum_{j=1}^n X_j , \sum_{k=1}^m Y_k \right)^2 \leq
mn \sum_{1 
\leq j \leq n, 1 \leq k \leq m} \tgm(X_i, Y_j)^2.
\end{equation}

\item
By operator concavity of the square root \cite[V.1.8, V.2.5]{Bhatia},
$X \mapsto \tgm(X,Y)$ is a concave function of $X$ (for fixed $Y$).
\end{enumerate}
\end{remarks}

\begin{lemma}
We can define a norm on the direct sum $H^\ell = H \oplus H \oplus
\cdots \oplus H$ of $\ell$ copies of $H$ by
\[
\|(\xi_1, \xi_2, \ldots, \xi_\ell)\|_{S1} = \trace \sqrt{ ( \langle
\xi_i, \xi_j \rangle )_{i,j=1}^\ell }
\]
\end{lemma}

\begin{proof}
Fix a unit vector $\eta \in H$ and consider the element $x \in
\mathcal{B}(H^\ell) \equiv M_\ell(\BH)$ given by
\[
x = \left(
\begin{array}{cccc}
\xi_1 \otimes \eta^* & \xi_2 \otimes \eta^* & \cdots & \xi_\ell
\otimes \eta^*\\[5mm]
\multicolumn{4}{c}{0}\\[5mm]
\end{array}
\right)
\]
This depends linearly on $(\xi_1, \xi_2, \ldots, \xi_\ell) \in
H^\ell$.
Computing $x^*x$ we find the block matrix $( \langle \xi_j, \xi_i
\rangle \eta \otimes \eta^* )_{i, j=1}^\ell$. Since $\eta \otimes
\eta^*$ is a self-adjoint projection, the square root of $x^*x$
has entries that are scalar multiples of the projection, the scalars
being the entries of the square root of $( \langle \xi_j, \xi_i 
\rangle )_{i, j=1}^\ell$. Since the projection has rank one,
it follows that the trace class norm of $x$
is the same as 
\[
\trace \sqrt{ ( \langle
\xi_i, \xi_j \rangle )_{i,j=1}^\ell }.
\]
Thus this expression gives a norm on $H^\ell$.
\end{proof}

\begin{example}
\rm
If $(\xi_1, \xi_2, \ldots, \xi_\ell) = (\lambda_1 \xi, \lambda_2 \xi,
\ldots, \lambda_\ell \xi)$ (for a unit vector $\xi \in H$ and scalars
$\lambda_j$), that is if the $\xi_j$ are linearly dependent, then
\[
\| (\xi_1, \xi_2, \ldots, \xi_\ell) \|_{S1} = \left( \sum_{j=1}^\ell
\|\xi_j\|^2 \right)^{1/2}.
\]

On the other hand if the $\xi_j$ are mutually orthogonal,
$\| (\xi_1, \xi_2, \ldots, \xi_\ell) \|_{S1} = \sum_{j=1}^\ell 
\|\xi_j\|$.
\end{example}

\begin{notation}
For $\mathbf{b} = [b_1, b_2, \ldots, b_\ell]^t \in \BH^\ell$, we may
regard $\mathbf{b}$ as an operator from $H$ to $H^\ell$ and then we
denote by $\|\mathbf{b}\|_{S1}$ the operator norm of $\mathbf{b}$
as an operator from $H$ to $(H^\ell, \| \cdot \|_{S1})$.
\end{notation}

\begin{theorem}
\label{ThmFormulaS1}
For $\mathbf{a} = [a_1, a_2, \ldots, a_\ell] \in \BH^\ell$,
$\mathbf{b} = [b_1, b_2, \ldots, b_\ell]^t \in \BH^\ell$
and $Tx = \sum_{j=1}^\ell a_j x b_j$,
we have
\[
\|T\| = \sup \left\{
\left\| \sqrt{Q(\mathbf{b}, \eta)^t} \mathbf{a}^* \right\|_{S1}
 : \eta \in H, \|\eta\| = 1
\right\}
\]
\end{theorem}

\begin{proof}
We use Theorem~\ref{thmformulaA}
and show that, for a fixed unit vector
$\eta \in H$
\[
\sup \left\{
\tgm \left( 
Q(\mathbf{a}^*, \xi), Q(\mathbf{b}, \eta)
\right)
: \xi \in H, \|\xi\| = 1
\right\}
=
\left\| \sqrt{Q(\mathbf{b}, \eta)^t} \mathbf{a}^* \right\|_{S1}
\]
It will be convenient to work with invertible $Q(\mathbf{b}, \eta)$
and it is sufficient to consider this case because of a perturbation
argument. For example, consider $\BH \subset \mathcal{B}(H^{\ell+1})$
with the $\ell+1$ copies of $H$ numbered $0$ to $\ell$ and $\BH$
included as the top left block in $\mathcal{B}(H^{\ell+1}) =
M_{\ell+1}(\BH)$. Let $E_{0j}$
denote the operator $E_{0j}(\xi_0, \xi_1, \ldots, \xi_\ell) = (0, 0,
\ldots, 0, \xi_0, 0, \ldots, 0)$ with the $\xi_0$ in position $j$. We
can replace $b_j$ by $b_j + \varepsilon E_{0j}$ for $\varepsilon> 0$ small.

From Theorem~\ref{thmformulaA}
we know that for $Sx = (Tx)p_2 = (Tx) (\eta
\otimes \eta^*)$
\[
\|S\| = \sup \left\{
\tgm \left( 
Q(\mathbf{a}^*, \xi), Q(\mathbf{b}, \eta)
\right)
: \xi \in H, \|\xi\| = 1 \right\}.
\]
We can rewrite $Sx$ using $\mathbf{a} \odot \mathbf{b}p_2 = (\mathbf{a} 
\alpha) \odot (\alpha^{-1} \mathbf{b}p_2)$ with $\alpha =
\sqrt{Q(\mathbf{b}, \eta)^t}$ to transform to the case where
$Q(\alpha^{-1}\mathbf{b}, \eta) = I_\ell$ is the $\ell \times \ell$ identity
matrix.

After this transformation,
\[
\|S\| = \sup \{ \trace \sqrt{Q((\mathbf{a} \alpha)^*, \xi)} : \xi \in H,
\|\xi\| = 1 \} = \| (\mathbf{a} \alpha)^* \|_{S1}
\]
and we have the result.
\end{proof}

\begin{lemma}
\label{LinIndepLemma}
If $\mathbf{b} = [b_1, b_2, \ldots, b_\ell]^t$ is an $\ell$-tuple of
elements of $\BH$ ($\ell \geq 1$), then $b_1, b_2, \ldots, b_\ell$ are
linearly independent if and only of there exist vectors $\xi_1,
\xi_2, \ldots, \xi_m \in H$ such that
\[
\sum_{k=1}^m Q(\mathbf{b}, \xi_k)
\]
is positive definite.
\end{lemma}

\begin{proof}
If the $b_j$ are linearly dependent, so that $\sum_{j=1}^\ell \lambda_j
b_j =0$ for some scalars $\lambda_j$ not all zero, then with $\lambda
= [\lambda_1, \lambda_2, \ldots, \lambda_\ell]$, for each $\xi
\in H$ we have
\[
\lambda Q(\mathbf{b}, \xi)^t 
\lambda^* =
\sum_{j,k=1}^\ell \lambda_j \bar{\lambda}_k \langle b_j \xi, b_k \xi
\rangle
= \left \| \sum_{j=1}^\ell \lambda_jb_j \xi \right\|^2 = 0
\]
and so each finite sum $\sum_{k=1}^m Q(\mathbf{b}, \xi_k)$ is
singular.

Conversely, if each finite sum is singular, choose a sum
$M=\sum_{k=1}^m Q(\mathbf{b}, \xi_k)$ of maximal rank among all such
sums and a nonzero vector
$\lambda \in \IC^\ell$ with $M^t\lambda^* =0$.
This $\lambda^*$
must be in the kernel of each $Q(\mathbf{b}, \xi)^t$ for $\xi \in H$
(as otherwise $Q(\mathbf{b}, \xi) + \sum_{k=1}^m Q(\mathbf{b},
\xi_k)$ would have larger rank).  From the above calculation we get
$\sum_{j=1}^\ell \lambda_jb_j \xi = 0$ for all $\xi$, hence
$\sum_{j=1}^\ell \lambda_jb_j =0$.
\end{proof}

We now characterise compactness of $T \in \El(\BH)$. Let
$H_1 = \{\xi \in H : \| \xi\| \leq 1\}$ denote the closed unit ball
of $H$ (compact in the weak topology) and 
let $\mathcal{K}(H)  $ denote the compact elements of $\BH$.

The equivalence (in the linearly independent case)
of the first and last conditions in the following
is known --- see \cite{FongSourour} or \cite[5.3.26]{AraMatBk} and
in the case $\ell=1$ see \cite{Vala}.

\begin{theorem}
\label{BHcompactness}
Let $\mathbf{a} = [a_1, a_2, \ldots, a_\ell] \in \BH^\ell$,
$\mathbf{b} = [b_1, b_2, \ldots, b_\ell]^t \in \BH^\ell$
and $Tx = \sum_{j=1}^\ell a_j x b_j$.
Then the following are
equivalent
\begin{enumerate}[(i)]
\item
\label{cpcti}
$T \colon \BH \to \BH$ is compact.
\item
\label{cpctii}
For $f_T \colon H_1 \times H_1 \to \IC$
given by
\[
f_T(\xi, \eta) = \tgm \left(
Q(\mathbf{a}^*, \xi), Q(\mathbf{b}, \eta)
\right),
\]
the function $f_T$ is continuous
in the product weak topology on $H_1 \times H_1$.
\item
\label{cpctiii}
The same function $f_T$ is continuous at each point of 
\[
(H_1 \times \{0\}) \cup
(\{0\} \times H_1) 
\]
in the product weak topology of $H_1 \times H_1$.
\end{enumerate}
Assuming that $a_j$ ($1 \leq j
\leq \ell$) are linearly independent and that $b_j$ ($1 \leq j 
\leq \ell$) are also linearly independent, then these are also
equivalent to
\begin{enumerate}
\item [(iv)]
$a_j, b_j \in \mathcal{K}(H)$
(for $1 \leq j \leq \ell$).
\end{enumerate}
\end{theorem}

\begin{proof} It is clear that (\ref{cpctii}) $\Rightarrow$
(\ref{cpctiii}). For the case of unit vectors $\xi, \eta$ we
have
\[
f_T (\xi, \eta) = \|T_{p_1, p_2}\|
\]
for $p_1 = \xi \otimes \xi^*$ and $p_2 = \eta \otimes \eta^*$,
and this is also the norm of the linear functional $x \mapsto \langle
(Tx)\eta, \xi \rangle$ by the proof of Lemma~\ref{firstlemma}.
We can see by homogeneity that for all $(\xi, \eta)$, $f_T(\xi, \eta)$
is the norm of the linear functional $x \mapsto \langle
(Tx)\eta, \xi \rangle$
and $f_T$ does not depend on the representation of $T$.
So it
is sufficient to establish the equivalence of the conditions under the
assumption that the $a_j$ ($1 \leq j \leq \ell$) and the $b_j$
($1 \leq j \leq \ell$) are linearly independent.

(\ref{cpcti}) $\Rightarrow$ (\ref{cpctii}):
As the
image $T(B_1)$ of the unit ball $B_1 = \{ x \in \mathcal{K}(H)
: \|x\| \leq 1 \}$
is relatively compact, for each $\varepsilon > 0$ we can find a finite
number of elements $x_1, x_2, \ldots, x_n \in T(B_1) \subset
\mathcal{K}(H)$ so that
\[
T(B_1) \subset \bigcup_{k=1}^n \{ y \in \mathcal{K}(H) : \|y - x_k\| <
\varepsilon \}.
\]
Denote by $\omega^T_{\xi, \eta}$ the linear functional on $\BH$ given
by $\omega^T_{\xi, \eta} (x) = \langle (Tx)\eta, \xi \rangle$ and
recall that $\|\omega^T_{\xi, \eta}\| = f_T(\xi, \eta)$.
The norm of $\omega^T_{\xi, \eta}$ is the same as the norm of its restriction
to $\mathcal{K}(H)$ (by
weak* continuity on $\BH$). Hence
\[
\max_{1 \leq k \leq n} | \langle x_k \eta, \xi \rangle |
\leq
f_T (\xi, \eta)
\leq
\left(\max_{1 \leq k \leq n} | \langle x_k \eta, \xi \rangle | \right) +
\varepsilon
\]
and so we have $F_\varepsilon \colon H_1 \times H_1
\to \IR$ given by $F_\varepsilon(\xi, \eta) = \max_{1 \leq k \leq
n} | \langle x_k \eta, \xi \rangle |$ with
\[
\sup_{(\xi, \eta) \in H_1 \times H_1} |f_T (\xi, \eta) -
F_\varepsilon(\xi, \eta)| < \varepsilon.
\]
For $x = \theta_1 \otimes \theta_2^* \in \BH$ of rank 1
($\theta_1, \theta_2 \in H$),
$(\xi, \eta) \mapsto
\langle x \eta, \xi \rangle = \langle \eta, \theta_2 \rangle \langle
\theta_1, \xi \rangle$ is clearly continuous on $H_1 \times H_1$. So then
is $(\xi, \eta) \mapsto
\langle x \eta, \xi \rangle$ when $x$ is of finite rank or when $x \in
\mathcal{K}(H)$ by approximation. Hence $F_\varepsilon$ is continuous
on $H_1 \times H_1$.
As a uniform limit of continuous functions, $f_T$ must be continuous.

(\ref{cpctiii}) $\Rightarrow$ (iv):
Using the linear independence assumption and Lemma~\ref{LinIndepLemma}
we can find $\eta_1, \eta_2, \ldots, \eta_n \in H$ so that
$\sum_{k=1}^n Q(\textbf{b}, \eta_k)$ is positive definite. We can
scale the $\eta_k$ so that $\eta_k \in H_1$ for $1 \leq k \leq n$.

Now, (weak) continuity of $f_T(\xi, \eta)$ at points $(0, \eta)$ implies
continuity of
\[
F_T(\xi, \eta) = \trace ( Q(\textbf{a}^*, \xi) Q(\textbf{b}, \eta) )
\]
at the same points because
\[
f_T(0, \eta) = F_T(0, \eta) = 0 \leq F_T(\xi, \eta) \leq \ell f_T(\xi,
\eta)^2.
\]
Thus we also have (weak) continuity at $\xi =0$ of 
$ \xi \mapsto \sum_{k=1}^n  F_T(\xi, \eta_k) $.
But as there is some $c > 0$ so that $\sum_{k=1}^n Q(\textbf{b},
\eta_k) > c I_\ell$, we have
\[
\trace  Q(\textbf{a}^*, \xi) \leq \frac{1}{c} \sum_{k=1}^n  F_T(\xi,
\eta_k) \to 0 \mbox{ as } \xi \to 0 \mbox{ weakly.}
\]
This means that $\sum_{j=1}^\ell \langle a_j a_j^* \xi, \xi \rangle
\to 0$ as $\xi \to 0$ weakly, or $\|a_j^* \xi \| \to 0$ as $\xi \to 0$
weakly. Hence each $a_j^*$ maps weakly null bounded sequences to norm
null sequences and $a_j^*$ 
 is compact. It follows that $a_j$ is compact. A similar argument
shows that each $b_j$ is compact.

(iv) $\Rightarrow$ (\ref{cpcti}): This follows from Vala's theorem.
\end{proof}

\section{General $C^*$-algebras}

We now extend some of 
the formulae to the case of $T \in \El(A)$
with $A$ an arbitrary $C^*$-algebra. For the remainder of
this section,
$Tx = \sum_{j=1}^\ell a_j x b_j$ for $x \in A$ with $a_j , b_j \in M(A)$.

By $P(A)$ we denote the set of pure states of $A$ and for $\phi \in P(A)$
let $\pi_\phi \colon A \to \mathcal{B}(H_{\pi_\phi})$ denote the associated
irreducible representation of $A$ (arrived at by the GNS method).
Observe that $\phi(x) = \langle \pi_\phi(x) \xi_\phi, \xi_\phi \rangle$
for $\xi_\phi$ the cyclic vector for the representation.
$\widehat{A}$ denotes the (unitary equivalence classes of) irreducible
representations of $A$ and we write $\pi_1 \sim \pi_2$ to indicate
unitary equivalence of representations. By $\phi_1 \sim \phi_2$ for
pure states, we mean equivalence of the associated representations.

Let
$T_\pi \colon \mathcal{B}(H_\pi) \to \mathcal{B}(H_\pi)$ be the
elementary operator 
\[
T_\pi(y) = \sum_{j=1}^\ell \pi(a_j) y \pi(b_j),
\]
using the fact that $\pi$ extends to $M(A)$. It is known that $\|T\| =
\sup_{\pi \in \widehat{A}} \|T_\pi\| $ (for example
\cite[5.3.12]{AraMatBk}).

For $\mathbf{b} = [b_1, b_2, \ldots, b_\ell]^t$ a column of elements
$b_j \in M(A)$  and $\phi \in P(A)$
we introduce the notation
\[
Q(\mathbf{b}, \phi) = ( \phi( b_i^* b_j ) )_{i, j =
1}^\ell
= \phi^{(\ell)} (Q(\mathbf{b}))
\in M_\ell(\IC).
\]
Here $\phi^{(\ell)} \colon M_\ell(A) \to M_\ell(\IC)$ is given by
$
\phi^{(\ell)} \left( (x_{ij})_{i,j=1}^\ell \right) =
(\phi(x_{ij}))_{i,j=1}^\ell.
$

\begin{theorem}
\label{CStarFormula}
For $T \in \El(A)$, $Tx = \sum_{j=1}^\ell a_j x b_j$,
$\mathbf{a} =[a_1, a_2, \ldots, a_\ell] \in M(A)^\ell$ a row and
$\mathbf{b} = [b_1, b_2, \ldots, b_\ell]^t$ a column of elements of
$M(A)$, we have
\begin{eqnarray*}
\|T\|
&=&
\sup \left\{
\tgm \left(
Q(\mathbf{a}^*, \phi_1),
Q(\mathbf{b}, \phi_2)
\right)
: \phi_1, \phi_2 \in P(A), \phi_1 \sim
\phi_2
\right\}
\\
&=& \sup \left\{
\left\| \sqrt{Q(\mathbf{b}, \phi)^t} \pi_\phi(\mathbf{a}^*) \right\|_{S1}
: \phi \in P(A)
\right\}
\end{eqnarray*}
(where $\pi_\phi(\mathbf{a}^*)$ means $[\pi_\phi(a_1^*),
\pi_\phi(a_2^*), \ldots, \pi_\phi(a_\ell^*) ]^t$).
\end{theorem}

\begin{proof}
This is immediate from the above remarks together with
Theorems~\ref{thmformulaA} and
\ref{ThmFormulaS1}.
\end{proof}

\begin{notation}
\rm
We now consider the maps $T^{(k)} \colon M_k(A) \to M_k(A) $ on the
space $M_k(A)$ of $k \times k$ matrices with entries in $A$, given
by $T^{(k)}( (x_{ij})_{i, j = 1}^\ell ) = (Tx_{ij})_{i, j = 1}^\ell$.
We use the canonical $C^*$ norms on $M_k(A)$ and the notation
$\|T\|_k = \left\| T^{(k)} \right\|$. The completely bounded norm
of $T$ is $\|T\|_{cb} = \sup_k \|T\|_k$.

Let $F_k(A)$ denote the set of factorial states considered in
\cite{ArchBatty}. They are those states $\phi$ of $A$ that are convex
combinations of at most $k$ unitarily equivalent pure states, or 
those where the commutant $\pi_\phi(A)'$ of the GNS representation is
a type $I_n$ factor with $ n \leq k$.
For a pure state $\psi$ on
$M_k(A)$ there is a factorial state $\phi \in F_k(A)$ with
\[
\psi \left(
\begin{array}{cccc}
x & 0 & \cdots & 0\\
0 & x &        & 0\\
  &   & \ddots & \\
0 & 0 & \cdots & x
\end{array}
\right) = \phi(x)
\]
for $x \in M(A)$. Conversely, if $\phi = \sum_{j=1}^k t_j \psi_j$ is a
convex combination of $\psi_j \in P(A)$ ($t_j \geq 0$, $\sum_j t_j
=1$) and we take the irreducible representation $\pi \colon A \to
\mathcal{B}(H_\pi)$ corresponding to
$\psi_1$, then there are unit vectors $\xi_j \in H$ so that $\psi_j(x)
= \langle x \xi_j, \xi_j \rangle$. The unit vector $(\sqrt{t_j}
\xi_j)_{j=1}^k \in H^k$ gives rise to a vector state $\psi$ for
$M_k(A)$ (acting on $H^k$ via $\pi^{(k)}$). This pure state $\psi$
will relate to $\phi$ as above.

From \cite[2.1 (iii)]{ArchBatty} we know that if $\phi$ is a proper convex
combination of states in $P(A)$, not all of which are equivalent, then
$\phi$ cannot be factorial. Relying on this, we can say that for
$\phi_1, \phi_2 \in F_k(A)$, $(\phi_1 + \phi_2)/2$ is factorial if and
only if the pure states in a convex combination making $\phi_1$ are
each
unitarily equivalent to those in a convex combination making $\phi_2$.
We write $\phi_1 \asymp \phi_2$ to mean that $(\phi_1 + \phi_2)/2$
is factorial.

From Theorem~\ref{CStarFormula}, we can deduce the following.
\end{notation}

\begin{corollary}
\label{CStarFormulak}
For $T \in \El(A)$, $Tx = \sum_{j=1}^\ell a_j x b_j$,
$
\mathbf{a} =[a_1, a_2, \ldots, a_\ell] \in M(A)^\ell
$
a row and
$\mathbf{b} = [b_1, b_2, \ldots, b_\ell]^t$ a column of elements of
$M(A)$, we have
\begin{eqnarray*}
\|T\|_k
&=&
\sup \left\{
\tgm \left(
Q(\mathbf{a}^*, \phi_1), Q(\mathbf{b}, \phi_2)
\right)
: \phi_1, \phi_2 \in F_k(A), \phi_1 \asymp
\phi_2
\right\}.
\end{eqnarray*}
\end{corollary}

\begin{corollary}
\label{CorollGrowthEst}
For $T \in \El(A)$,
$Tx =
\sum_{j=1}^\ell a_j x b_j$, we have
\begin{eqnarray*}
\|T\|_k & \leq & \max(k, \sqrt{\ell}) \|T\|\\
\|T\|_{cb} & \leq & \sqrt{\ell} \|T\|
\end{eqnarray*}
\end{corollary}

\begin{proof}
It is well-known that (for any operator $T \colon A \to A$) $\|T\|_k
\leq k \|T\|$. (See \cite[Exercise 3.10(ii)]{Paulsen}.)

In \cite{TimoneyNE}, it was shown that $\|T\|_{cb} =
\|T\|_\ell$ and so we could deduce $\|T\|_k \leq \ell \|T\|$ for all
$k$, but we seek the improved bound involving $\sqrt{\ell}$.

With $X = Q(\mathbf{a}^*, \phi_1)$, $Y = Q(\mathbf{b}, \phi_2)$,
$\phi_1, \phi_2 \in P(A)$, $\phi_1 \sim \phi_2$, 
Theorem~\ref{CStarFormula} tells us
\[
\tgm(X, Y) = \sum_{j=1}^\ell \sqrt{\lambda_j(XY)} \leq \|T\|
\Rightarrow \sum_{j=1}^\ell \lambda_j(XY) = \trace (XY) \leq \|T\|^2
\]
This latter is a convex condition on $X$ and $Y$. Hence it remains
true on replacing $\phi_1, \phi_2 \in P(A)$ by convex combinations.
Hence for
$\phi_1, \phi_2 \in F_k(A)$ with $\phi_1 \asymp \phi_2$,
\[
\trace Q(\mathbf{a}^*, \phi_1)  Q(\mathbf{b}, \phi_2) \leq \|T\|^2
\Rightarrow \tgm( Q(\mathbf{a}^*, \phi_1), Q(\mathbf{b}, \phi_2)) \leq
\sqrt{\ell} \|T\|
\]
(by the Cauchy Schwarz inequality for $\sum_j \sqrt{\lambda_j}$). Thus
Corollary~\ref{CStarFormulak} implies $\|T\|_k \leq \sqrt{\ell}
\|T\|$.
\end{proof}

\begin{example}
\rm
\begin{enumerate}[(i)]
\item
The well-known
transpose example $T \colon M_n \to M_n$, $Tx = x^t$ has
$\|T\| = 1$ while $\|T\|_{cb} = \|T\|_n = n$. As an elementary
operator $Tx = \sum_{i,j=1}^n e_{ij}x e_{ij}$ (with $e_{ij}$ the
matrix with 1 in the $(i, j)$ place and zeros elsewhere) and has $\ell
= n^2$. Thus the familiar estimate $\|T\|_k \leq k \|T\|$ cannot be
improved in general.

For $Sx = e_{11}x^t = \sum_{j=1}^n e_{1j} x e_{1j}$ the first row of
the transpose, one can check that $\|S\|_n = \sqrt{n}$ (and $\|S\| =
1$) so that the $\sqrt{\ell}$ bound is also optimal in a sense.

For this $S$, it is true that $\|S\|_k = \sqrt{k}$
for $1 \leq k \leq n$, while for the transpose example
$\|T\|_k = k$ ($1 \leq k \leq n$).

\item
A natural question would be to describe the functions of $k$ that can
arise as $f(k) = \|T\|_k$ for $T \in \El(A)$.
This would mean finding
further refinements of the estimates we have to pin
down the possibilities for $f(k)$.
For example $f(k_1k_2)
\leq k_1 f(k_2)$ follows because $T^{(k_1k_2)} = (T^{(k_2)})^{(k_1)}$.
It follows that, for $n \leq k$, $f(k+n) \leq f(2k) \leq 2 f(k)$.

If $T_1 \in \El(A_1)$ and $T_2 \in \El(A_2)$ then $T_1 \oplus T_2 \in
\El(A_1 \oplus A_2)$ satisfies $\|T_1 \oplus T_2\|_k = \max(
\|T_1\|_k, \|T_2\|_k)$. Hence we can, for example, combine the
identity with the
transpose example. Take $T_1x = x$ and
$T_2x = x^t/m$ on $M_n(\IC)$ (where $1 < m < n$) to produce $T = T_1
\oplus T_2$ where $f(k) = \|T\|_k = 1$ for $1 \leq k \leq m$, while
$f(k) = k/m$ for $m+1 \leq k \leq n$. In examples of this type
$f(k+1)/f(k)$ can be 1 and $(k+1)/k$ in different intervals.

We can use (\ref{quasiconvexity}) to show that
\begin{equation}
\label{StepGrowthEstim}
\|T\|_{k+1} \leq \left( 1 + \frac{2 \sqrt{k}}{k+1} \right) \|T\|_k.
\end{equation}
holds for general $T \in \El(A)$ and $k = 1, 2, \ldots$. Assume $\|T\|_k
\leq 1$ and consider a pair of factorial states $\phi_i = \sum_{r=1}^{k+1}
t_{ir} \phi_{ir}$ where $\sum_r t_{ir} = 1$, $t_{ir} \geq 0$,
$\phi_{ir}$ are pure states that are all unitarily equivalent ($i =1,
2$, $1 \leq r \leq k+1$). For at least one $r$, $t_{ir} \leq 1/(k+1)$
and we assume $t_{i1} \leq 1/(k+1)$ for $i = 1, 2$. Then write
$\phi_i = t_{i1} \phi_{i1} + (1-t_{i1}) \psi_i$ for 
$\psi_i \in F_k(A)$ ($i = 1, 2$). Let $X = Q(\mathbf{a}^*, \phi_1)$, $Y =
Q(\mathbf{b}, \phi_2)$, $X_1 = Q(\mathbf{a}^*, \phi_{1i})$, $X_2 =
Q(\mathbf{a}^*, \psi_1)$, $Y_1 = Q(\mathbf{b}, \phi_{21})$, $Y_2 =
Q(\mathbf{b}, \psi_2)$ so that $X = t_{11} X_1 + (1- t_{11})X_2$ and
$Y = t_{21}Y_1 + (1- t_{21}) Y_2$. Using (\ref{quasiconvexity})
and $\tgm(X_i, Y_\alpha) \leq 1$ for $i, \alpha = 1, 2$ we get
\begin{eqnarray*}
\tgm(X, Y) &\leq&
(\sqrt{t_{11}} + \sqrt{1 - t_{11}})
(\sqrt{t_{21}} + \sqrt{1 - t_{21}})\\
&\leq& \left(\sqrt{\frac{1}{k+1}} + \sqrt{\frac{k}{k+1}} \right)^2.
\end{eqnarray*}
By Corollary~\ref{CStarFormulak} we get (\ref{StepGrowthEstim}).
\end{enumerate}
\end{example}

The following result is shown in \cite{AMS} but the proof
relies in an essential way on \cite[Theorem 3.1]{Magajna} (concerning
the case of prime $C^*$-algebras with zero socle).

\begin{corollary}
\label{antiliminalcoroll}
If $A$ is an antiliminal $C^*$-algebra and $T \in \El(A)$,
then $\|T\| = \|T\|_{cb}$.
\end{corollary}

\begin{proof}
It is known that for antiliminal
$C^*$-algebras (or more generally for `antiliminal by abelian' ones
\cite{ArchBatty})
every factorial state is a weak*-limit of pure states.
To deduce the Corollary directly from
Corollary~\ref{CStarFormulak}, we would need to know further
that for
$\phi_1, \phi_2 \in F_k(A)$ with $\phi_1 \asymp \phi_2$, we can find
a net of pairs of pure states
$\phi_{1, \alpha}, \phi_{2, \alpha} \in P(A)$ with 
$\lim_\alpha \phi_{1, \alpha} = \phi_1$,
$\lim_\alpha \phi_{2, \alpha} = \phi_2$ and $\phi_{1, \alpha} \sim
\phi_{2, \alpha}$ for each $\alpha$.

However, consideration of the proof in \cite[11.2.3]{DixCstar}
reveals that this further fact is true. From the Kadison transitivity
theorem,
we can choose a pure state $\phi$, unitary $u_{ij}
\in \tilde{A} =$ the unitisation of $A$ and $0 \leq t_{ij} \leq 1$
($i = 1,2 $, $1 \leq j \leq k$) so that
\[
\phi_i(x) = \sum_{j=1}^k t_{ij} \phi(u_{ij}^* x u_{ij}) ,
\sum_{j=1}^k t_{ij} = 1 \quad (i = 1, 2).
\]
The state $\phi$ can then be approximated by a net $\psi_\alpha$ of
states where there exists an irreducible representation
$\pi_\alpha \colon A \to \mathcal{B}(H_{\pi_\alpha})$ with $\psi_\alpha$
vanishing on the inverse image
$\pi_\alpha^{-1}(\mathcal{K}(H_{\pi_\alpha}))$ of
the compacts (\cite[11.2.2]{DixCstar}).
By \cite[11.2.1]{DixCstar}, the two states
\[
\psi_{i\alpha} (x) 
=
\sum_{j=1}^k t_{ij} \psi_\alpha (u_{ij}^* x u_{ij})
\quad (i = 1, 2)
\]
can be approximated by vector states on $\mathcal{B}(H_{\pi_\alpha})$
composed with $\pi_\alpha$, hence by equivalent pure states of $A$.
\end{proof}

\section{Compact elementary operators}

We now characterise compact elementary
operators on a $C^*$-algebra $A$. For the case $A = \BH$,
Theorem~\ref{BHcompactness} extends known results somewhat.  For
prime $C^*$-algebras $A$, \cite[5.3.26]{AraMatBk} gives a
characterisation but a similar result for general $A$ seems not to
have been established up to now. In addition, we characterise weakly
compact elements of $\El(\BH)$ in a similar way to
Theorem~\ref{BHcompactness}.

Compactness of elements $a$ of a $C^*$-algebra $A$ (in terms of $a$
being in the closure of the socle, or weak compactness of the
operators of left or right multiplication by $a$, or compactness of the
elementary operator $x \mapsto axa$) has been studied by several
authors and references can be found in \cite[p.~36]{AraMatBk}.
The set $\mathcal{K}(A)$ of compact elements of $A$ is a closed ideal in
$A$.

We introduce the notion of the $\mathcal{K}(A)$ topology on the space
$E(A) = \{ \psi : \psi =\lambda \phi , \phi \in P(A), 0 \leq \lambda
\leq 1 \}$ of multiples of pure states (a subset of the dual of $A$).
This means the topology of pointwise convergence on $\mathcal{K}(A)$,
so that a net $(\psi_\alpha)_\alpha$ converges to $\psi \in E(A)$ if
and only if $\lim_\alpha \psi_\alpha(x) = \psi(x)$ for each $x \in
\mathcal{K}(A)$. If $\mathcal{K}(A)$ is small (zero for example) then
this will be a very coarse topology on $E(A)$ (even the
trivial topology if $\mathcal{K}(A)=0$).
By $R(A)$ we denote the subset of the product $E(A)
\times E(A)$ consisting of pairs $(\lambda_1 \phi_1, \lambda_2
\phi_2)$ with $\phi_1 \sim \phi_2$ unitarily equivalent pure states
(and $0 \leq \lambda_1, \lambda_2 \leq 1$).

When $A = \mathcal{K}(H)$, there is a surjection $\mu \colon H_1 \to
E(A)$ given by $\mu(\xi) = \omega_\xi$.
This is continuous
(from the weak topology on $H_1$ to the $\mathcal{K}(A)$
topology on $E(A)$).
Neighbourhoods
of $\omega_\xi \in E(A)$ contain finite intersections of
neighbourhoods of the form $N= \{\omega_\eta : |\omega_\eta(x) -
\omega_\xi(x)| < 1\}$ with $x = \theta_1 \otimes \theta_2^*$ of rank
one ($\theta_1, \theta_2 \in H$). Since
$\mu^{-1}(N) = \{ \eta \in
H_1: |<\eta, \theta_2><\theta_1, \eta> -
<\xi, \theta_2><\theta_1, \xi>| < 1 \}$,
continuity of $\mu$ is clear.
$\mu$ is also surjective and the inverse image of $\omega_\xi$
consists of multiples $\zeta \xi$ with $\zeta \in \mathbb{T} = \{
\zeta \in \IC: |\zeta| = 1\}$.
As the quotient space of $H_1$ by the action of $\mathbb{T}$ is compact
in the quotient topology and $E(A)$ is Hausdorff,
it follows that $\mu$ induces a
homeomorphism from the quotient space $H_1/\mathbb{T}$
to $E(A) = E(\mathcal{K}(H))$.

In the case $A = \BH$ (with $H$ infinite dimensional) the elements of
$E(\BH)$ not in $E(\mathcal{K}(H))$ are the nonzero multiples $\lambda
\phi$ of pure states $\phi$ of $\BH$ vanishing on $\mathcal{K}(H)$,
and so they are
not separated from $0$ in the $\mathcal{K}(\BH)$ topology of $E(\BH)$.

For the case of a general $C^*$-algebra $A$, we can also analyse the
topology on $E(A)$. Pure states $P(A)$ of $A$ are either pure states of
$\mathcal{K}(A)$ or vanish on $\mathcal{K}(A)$, so that $P(A) =
P(\mathcal{K}(A)) \cup P(A/\mathcal{K}(A))$ is a disjoint union.
Note that multiples of states in $P(A/\mathcal{K}(A))$
are in the closure of 0 in the $\mathcal{K}(A)$
topology on $E(A)$. Similarly  irreducible representations
$\pi \colon A
\to \mathcal{B}(H_\pi)$ of $A$ are either irreducible when restricted
to $\mathcal{K}(A)$ or may be regarded as irreducible representations
of the quotient. We will use the notation $\widehat{A}$ in a 
\emph{slightly
ambiguous way} as also standing for a specific set of representatives
of the unitary equivalence classes of irreducibles. We have then
$\widehat{A} = \widehat{\mathcal{K}}(A) \cup
(A/\mathcal{K}(A)) \hatted$ (where we use $\widehat{\mathcal{K}}(A)$
for $(\mathcal{K}(A)) \hatted$).
For $\pi \in  \widehat{\mathcal{K}}(A)$ we have $\pi(\mathcal{K}(A)) =
\mathcal{K}(H_\pi)$.

It is known that $\mathcal{K}(A)$ is a direct sum of algebras
$\mathcal{K}(H_i)$ of compact operators on Hilbert spaces $H_i$
(for $i \in I =$ some index set). This is known because $\mathcal{K}(A)$
is a $C^*$-algebra of compact elements and we can apply \cite[Theorem
8.2]{JCAlexander} (or a fact now known to be equivalent shown in
\cite[Theorem 8.3]{KaplanskyNormed}). It follows that the $H_i$ are
the $H_{\pi}$ with $\pi \in \widehat{\mathcal{K}}(A)$. Also for $\lambda \phi
\in E(A)$ with $\lambda \neq 0$ and
where $\phi \in P(\mathcal{K}(A))$, $\pi_\phi$ is (unitarily
equivalent to) some $\pi \in \widehat{\mathcal{K}}(A)$. Hence there is a nonzero
vector $\xi$ in the unit ball of $H_\pi$ so that $\lambda\phi(x) =
\omega_\xi(\pi(x))$. By considering $\mathcal{K}(H_\pi)$ as contained
in the direct sum $\mathcal{K}(A)$ we see that there exist neighbourhoods of
$\lambda\phi$ in $E(A)$ that consist entirely of functionals
$\psi \in E(\mathcal{K}(H_\pi)) \setminus \{0\}$.
So the
complement of the closure of 0 in
$E(A)$
is homeomorphic to the disjoint union of $E(\mathcal{K}(H_\pi)) \setminus
\{0\}$ (for $\pi \in \widehat{\mathcal{K}}(A)$).
Moreover
each $E(\mathcal{K}(H_\pi)) \setminus
\{0\}$ is homeomorphic to $((H_\pi)_1 \setminus \{0\})/\mathbb{T}$
(for $\pi \in \widehat{\mathcal{K}}(A)$).
On the other hand neighbourhoods of $0$ in $E(A)$ contain all but
finitely many of $E(\mathcal{K}(H_\pi))$ for $\pi \in \widehat{\mathcal{K}}(A)$
(and intersect the remaining $E(\mathcal{K}(H_\pi))$ in neighbourhoods
of 0). This is the case because a basic neighbourhood of $\psi \in
E(A)$ is an intersection of neighbourhoods of the form $N_{\psi, x} =
\{\psi_1 \in
E(A) : |\psi_1(x) - \psi(x)| < 1 \}$ with $x \in \mathcal{K}(A)$. Any
such $x$ can be approximated in norm by $x_0 \in \bigoplus_{j=1}^n
\mathcal{K}(H_{\pi_j})$,  a finitely supported element of the direct
sum making up $\mathcal{K}(A)$. Choosing $x_0$ so that $\|x - x_0\| <
1/4$ we get $N _{\psi, 2x_0} \subset N _{\psi,  x}$ and the smaller
neighbourhood places a restriction only on finitely many
$E(\mathcal{K}(H_\pi))$.

For $(\phi, \psi) \in R(A)$,
if $\phi$ is outside the closure $\overline{\{0\}}$
of 0 in the $\mathcal{K}(A)$ topology on
$E(A)$ then $\psi$ is either 0 or also
outside $\overline{\{0\}}$.  We can say that
$R(A)$ is the union of $R(\mathcal{K}(A))$ and a set contained in
$\overline{\{0\}} \times \overline{\{0\}}$.

\begin{theorem}
\label{CompactnessThm}
Let $A$ be a $C^*$-algebra and $T \in \El(A)$.
Then the following are equivalent for $T$
\begin{enumerate}[(i)]
\item
\label{Acpcti}
$T$ is compact
\item
\label{Acpctii}
If $Tx = \sum_{j=1}^\ell a_j x b_j$ for $a_j , b_j \in M(A)$, and $f_T
\colon R(A) \to \IR$ is defined as
\[
f_T(\phi_1, \phi_2) = \tgm \left(
Q(\mathbf{a}^*, \phi_1), Q(\mathbf{b}, \phi_2)
\right),
\]
the function $f_T$ is continuous on $R(A)$.

\item
\label{Acpctiii}
$T$ can be expressed
as $Tx = \sum_{j=1}^\ell a_j x b_j$ for $a_j$ and $b_j$ compact
elements of $A$ ($1 \leq j \leq \ell$).
\end{enumerate}
\end{theorem}

\begin{proof}

(\ref{Acpcti}) 
$\Rightarrow$
(\ref{Acpctii}):
Let $Tx = \sum_{j=1}^\ell a_j x b_j$.

Compactness of $T$ implies compactness of its double transpose (again
$x \mapsto \sum_{j=1}^\ell a_j x b_j$) and of its restriction $T_{\pi_a}$
to the
atomic part $\prod_{\pi \in \widehat{A}} \mathcal{B}(H_\pi)$ of the double
dual of $A$ (where the same formula $x \mapsto \sum_{j=1}^\ell a_j x
b_j$ holds modulo identifying elements of $M(A)$ with their images
under the reduced atomic representation $\pi_a = \bigoplus_{\pi \in \widehat{A}}
\pi$).
Let $T_\pi \in \El(\mathcal{B}(H_\pi))$ be
the elementary operator arising from $T$ and $\pi$. As a restriction
of $T_{\pi_a}$, $T_\pi$ is compact.  Moreover, for each $\varepsilon >
0$, we have $\{ \pi \in \widehat{A} : \|T_\pi\| > \varepsilon \}$
finite. This is a consequence of \cite[5.3.17]{AraMatBk}, but we offer
a self-contained argument.

If not we could find a sequence of unit norm elements
$x_n \in \mathcal{B}(H_{\pi_n})$ with $\pi_n \in \widehat{A}$ all distinct
and $\|T_{\pi_n}(x_n)\| > \varepsilon$. Inside $\prod_{\pi \in \widehat{A}}
\mathcal{B}(H_\pi)$ we have then an infinite dimensional space spanned
by the $x_n$ on which $T$ restricts to a topological linear
isomorphism, contradicting compactness of $T$.

As each $T_\pi$  ($\pi \in \widehat{A}$) is
compact, we have $T_\pi(\mathcal{B}(H_\pi)) \subset
\mathcal{K}(H_\pi)$ by Theorem~\ref{BHcompactness}. Putting these
facts together, we get $T(A) \subset \mathcal{K}(A)$ (via
\cite[1.2.30(g)]{AraMatBk}).

For $x \in A$ and $(\phi, \psi) \in R(A)$ we define
$|x(\phi, \psi)|$ as follows.
Choose $\pi \in
\widehat{A}$
and $\xi, \eta \in (H_\pi)_1$ with $\phi(x) = \omega_\xi(\pi(x))$,
$\psi(x) = \omega_\eta(\pi(x))$, and set $|x(\phi, \psi)| = |\langle \pi(x)
\eta, \xi \rangle|$.

Then continuity of $f_T$ can be established using $\|T\| = \sup \{ |y(\phi,
\psi)|: y \in T(A_1) , (\phi, \psi) \in R(A) \}$ (where $A_1$ is the
unit ball of $A$) and an argument
similar to the one used in the proof of (\ref{cpcti}) $\Rightarrow$
(\ref{cpctii}) of Theorem~\ref{BHcompactness}. Choose $y_1, y_2,
\ldots, y_n \in T(B_1) \subset \mathcal{K}(A)$ so that each $y \in
T(B_1)$ has $\|y - y_k\| < \varepsilon$ for some $1 \leq k \leq n$.
$f_T(\phi, \psi)$ is approximately $\max_k |y_k(\phi, \psi)|$.

To finish, we need to know that for $z \in \mathcal{K}(A)$ the map
\[
(\phi, \psi) \mapsto |z(\phi, \psi)| \colon R(A) \to \IR
\]
is
continuous. We can approximate $z$ in norm by a finitely supported
element
$z' = \sum_{j=1}^n z_j \in \bigoplus_{j=1}^n \mathcal{K}(H_{\pi_j})$
in the direct sum making up $\mathcal{K}(A)$. Then $|z'(\phi, \psi)| =
\sum_{j=1}^n |z_j(\phi, \psi)|$. The continuity of $(\phi, \psi)
\mapsto |z_j(\phi, \psi)|$ can be shown by
noting first that the map is zero on those parts of $R(A)$
coming from pairs of pure states with GNS representation different from
$\pi_j$. Consider the weakest topology $\tau$ on $R(A)$ so that the maps
$(\phi, \psi) \mapsto (\phi(z'), \psi(z'')) \colon R(A) \to \IC^2$ are
continuous for $z', z'' \in \mathcal{K}(H_{\pi_j})$. Note that $\tau$
is weaker than the usual topology on $R(A)$.
This space $(R(A), \tau)$ is
$R(\mathcal{K}(H_{\pi_j}))$ plus points that cannot be separated from
$(0, 0)$, and it becomes $R(\mathcal{K}(H_{\pi_j})) =
E(\mathcal{K}(H_{\pi_j})) \times E(\mathcal{K}(H_{\pi_j}))$ if we
identify these points with $(0, 0)$. We saw in the proof of
Theorem~\ref{BHcompactness} that $(\eta, \xi) \mapsto \langle  z_j
\eta, \xi \rangle$ is continuous on $(H_\pi)_1 \times (H_\pi)_1$ and
it follows then that $(\eta, \xi)
\mapsto | \langle  z_j
\eta, \xi \rangle |$ is continuous on
$(H_\pi)_1/\mathbb{T} \times (H_\pi)_1/\mathbb{T}$. Hence
$(\omega_\xi, \omega_\eta) \mapsto | \langle  z_j
\eta, \xi \rangle |$ is continuous on
$E(\mathcal{K}(H_{\pi_j})) \times E(\mathcal{K}(H_{\pi_j}))$
and $(\phi, \psi)
\mapsto |z_j(\phi, \psi)|$ is continuous on $(R(A), \tau)$, thus
on $R(A)$.

(\ref{Acpctii}) 
$\Rightarrow$
(\ref{Acpctiii}):
Let $T \in \El(A)$ and write $Tx = \sum_{j=1}^\ell a_{0j} x
b_{0j}$ for some $a_{0j}, b_{0j} \in M(A)$. Our aim is to show that there is an
alternative way to write $T$ with $a_{0j}$ and $b_{0j}$ replaced by compact
elements of $A$.

First, note that continuity of $f_T$ implies that $f_T(\phi, \psi)
= 0$ when $(\phi, \psi ) \in R(A/\mathcal{K}(A))$.
Since $\|T\|_\pi = \sup \{ f_T(\phi,
\psi) : (\phi, \psi) \in R(\mathcal{B}(H_\pi)) \}$, it follows that
$T_\pi = 0$ for $\pi \in (A/\mathcal{K}(A)) \hatted $.
Also note that
for each $\varepsilon > 0$, we have $\{ \pi \in \widehat{\mathcal{K}(A)}
: \|T_\pi\| > \varepsilon \}$ finite. This follows from continuity
of $f_T$ at $(0,0)$ and earlier remarks about neighbourhoods of 0 in
$E(A)$.

Let $S$ denote the restriction of $T$ to
$\mathcal{K}(A)$ and observe that $S \in \El(\mathcal{K}(A))$ and $S$
has a representation as a sum of $\ell$ terms.

The function $f_S$ is the restriction of $f_T$ to $R(\mathcal{K}(A))$
and the function $f_{S_\pi}$ of Theorem~\ref{BHcompactness}
on $(H_\pi)_1 \times (H_\pi)_1$ is
$f_{S_\pi}(\xi, \eta) = f_T(\omega_\xi \circ \pi, \omega_\eta \circ
\pi)$. Therefore
$f_{S_\pi}$ is
continuous by the relationships between the topologies concerned.
Thus each $S_\pi$ is compact for $\pi \in \widehat{\mathcal{K}}(A)$.

By Theorem~\ref{BHcompactness} or \cite[5.3.26]{AraMatBk}, each $S_\pi$
can be represented as $S_\pi(x) = \sum_{j=1}^\ell a_{\pi j} x b_{\pi j}$ with
$a_{\pi j}, b_{\pi j} \in \mathcal{K}(H_\pi)$. We can moreover arrange that
\[
\|S_\pi\|_{cb} = \left\| \sum_{j=1}^\ell a_{\pi j} \otimes b_{\pi j} \right\|_h
= \left\| \sum_{j=1}^\ell a_{\pi j} a_{\pi j}^* \right\|
= \left\| \sum_{j=1}^\ell b_{\pi j}^* b_{\pi j} \right\|
\]
using the Haagerup's theorem \cite[5.4.7]{AraMatBk} and the fact
that the infimum defining the Haagerup norm of a tensor can be
realised without increasing the length of its representation or the
span of either $\{a_{\pi j} : 1 \leq j \leq \ell\}$ or $\{b_{\pi j} : 1
\leq j \leq \ell\}$ \cite[Proposition 9.2.6]{EffrosRuan}.  By
Corollary~\ref{CorollGrowthEst}, $\|S_\pi \|_{cb} \leq \sqrt{\ell}
\|S_\pi \|$ and so $\max( \|a_{\pi j}\|^2, \|b_{\pi j}\|^2) \leq \sqrt{\ell}
\|S_\pi \|$. It follows that we can define $a_j, b_j \in \bigoplus_{\pi 
\in \widehat{\mathcal{K}}(A)} \mathcal{K}(H_\pi ) = \mathcal{K}(A)$
by $a_j = (a_{\pi j})_{\pi  \in
\widehat{\mathcal{K}}(A)}$ and $b_j = (b_{\pi j})_{\pi  \in \widehat{\mathcal{K}}(A)}$. Then we have $Sx = \sum_{j=1}^\ell
a_j x b_j$ for $x \in \mathcal{K}(A)$ since $\pi(Sx) = S_\pi(\pi(x))
= \sum_{j=1}^\ell a_{\pi j} \pi(x) b_{\pi j} = \pi \left(
\sum_{j=1}^\ell
a_j x b_j \right)$ for $\pi \in \widehat{\mathcal{K}}(A)$.

Let $T_1 \in \El(A)$ be $T_1 x = \sum_{j=1}^\ell a_j x b_j$. We have $T_1 x =
Sx = Tx$ for $x \in \mathcal{K}(A)$.

Finally we show $T_1 = T$ to complete the proof. Let $x \in A$.
For any
$\pi \in (A/\mathcal{K}(A)) \hatted $ we have
$\pi(Tx) = 0 = \pi(T_1x)$.
For $\pi \in \widehat{\mathcal{K}}(A)$, $\pi(\mathcal{K}(A)) =
\mathcal{K}(H_\pi)$. Since $(T_1)_\pi = T_\pi = S_\pi$ on
$\mathcal{K}(H_\pi)$, a density argument shows $(T_1)_\pi = T_\pi$.
Thus $\pi(Tx) = T_\pi(\pi(x))
= (T_1)_\pi(\pi(x)) = \pi(T_1x)$ for $\pi \in
\widehat{\mathcal{K}}(A)$ also.

(\ref{Acpctiii}) 
$\Rightarrow$
(\ref{Acpcti}):
If $a_j$ and $b_j$ are each compact ($1 \leq j \leq \ell$) then $T$ is
compact by Vala's theorem (as remarked in proof of
\cite[5.3.26]{AraMatBk}).
\end{proof}

We now consider weak compactness of elementary operators on $\BH$ (for
$H$ infinite dimensional).
Our arguments will involve the matrix valued essential numerical
ranges $Q(\mathbf{a}^*, \phi)$ and $Q(\mathbf{b}, \phi)$ where $\phi$
is a state of the Calkin algebra $\BH/\mathcal{K}(H)$. A fact we
rely upon is that all states $\phi$ of the Calkin algebra (or states
of $\BH$ vanishing on $\mathcal{K}(H)$) are
weak*-limits of vector states $\omega_\xi(x) = \langle x \xi, \xi
\rangle$ (for unit vectors $\xi \in H$) \cite[11.2.1]{DixCstar}.
That is, for any state $\phi$
vanishing on $\mathcal{K}(H)$
there is a net $\omega_{\xi_\alpha}$ so that $\phi(x) = \lim_\alpha
\omega_{\xi_\alpha}(x)$ for all $ x \in \BH$. Taking $x$ to be a rank
one operator $\eta^* \otimes \eta$ we see that $0 = \phi(\eta^*
\otimes \eta) = \lim_\alpha |\langle \xi_\alpha, \eta \rangle|^2$ and so 
$\xi_\alpha \to 0$ weakly in $H$.

\begin{lemma}
\label{EssLinIndepLemma}
If $\mathbf{b} = [b_1, b_2, \ldots, b_\ell]^t$ is an $\ell$-tuple of
elements of $\BH$ ($\ell \geq 1$), then $b_1, b_2, \ldots, b_\ell$ are
linearly independent modulo $\mathcal{K}(H)$
if and only if there exists a weakly null net $(\xi_\alpha)_\alpha$ of
unit vectors in $H$
such that
\[
\lim_\alpha Q(\mathbf{b}, \xi_\alpha)
\]
is positive definite.
\end{lemma}

\begin{proof}
We assume that the $\ell$-tuple of elements
$b_j + \mathcal{K}(H)$ ($1 \leq j \leq \ell$)
of the Calkin algebra are linearly independent.
We apply Lemma~\ref{LinIndepLemma}, by
temporarily representing $\BH/\mathcal{K}(H)$
as an algebra of operators on some Hilbert space. Normalising the
vectors we obtain from Lemma~\ref{LinIndepLemma} so that the sum of
the squares of their norms is 1, we see that we have
a state $\phi$ of $\BH/\mathcal{K}(H)$ so that $Q(\mathbf{b}, \phi)$
is positive definite.

Applying the above remarks about states of the Calkin algebra we
deduce the existence of the net $(\xi_\alpha)_\alpha$ as stated.

For the converse, given a weakly null net $(\xi_\alpha)_\alpha$  of
unit vectors as in the statement,
we can pass to a subnet and assume that the vector
states $\omega_{\xi_\alpha}$ converge weak* to some state $\phi$ of
$\BH$. One can easily see that $\phi$ vanishes on rank one operators
and so on $\mathcal{K}(H)$. Note that $Q(\mathbf{b}, \phi) =
\lim_\alpha Q(\mathbf{b}, \omega_{\xi_\alpha})$ is positive
definite. If a linear combination $\sum_{j=1}^\ell \lambda_j b_j \in
\mathcal{K}(H)$, then for $\lambda = [\lambda_1, \lambda_2, \ldots,
\lambda_\ell]$ we can compute
\[
0 = \phi \left( \left( \sum_{j=1}^\ell \lambda_j b_j \right)^*
\left( \sum_{j=1}^\ell \lambda_j b_j \right) \right) =
\lambda^* Q(\mathbf{b}, \phi)^t \lambda =0
\]
and so $\lambda = 0$.
\end{proof}

\begin{theorem}
\label{BHweakcompactness}
Let $\mathbf{a} = [a_1, a_2, \ldots, a_\ell] \in \BH^\ell$,
$\mathbf{b} = [b_1, b_2, \ldots, b_\ell]^t \in \BH^\ell$
and $Tx = \sum_{j=1}^\ell a_j x b_j$.
Then the following are
equivalent
\begin{enumerate}[(i)]
\item
\label{wcpcti}
$T \colon \BH \to \BH$ is weakly compact.
\item
\label{wcpctii}
$T \colon \mathcal{K}(H) \to \mathcal{K}(H)$ is weakly compact.
\item
\label{wcpctiii}
For $f_T \colon H_1 \times H_1 \to \IC$
given by
\[
f_T(\xi, \eta) = \tgm \left(
Q(\mathbf{a}^*, \xi), Q(\mathbf{b}, \eta)
\right),
\]
the function $f_T$ is continuous at $(0,0)$
(in the product weak topology on $H_1 \times H_1$, with $H_1$ the
closed unit ball of $H$ as before).
\item
\label{wcpctiv}
The function $g_T(\xi) = f_T(\xi, \xi)$ on $H_1$ is continuous at 0.
\item
\label{wcpctv}
There exist $c_1, c_2 , \ldots, c_\ell \in \mathcal{K}(H)$, $d_1,
d_2, \ldots , d_\ell \in \BH$ and $0 \leq m \leq \ell$ so that
\[
Tx = \sum_{j=1}^m c_j x d_j + \sum_{j=m+1}^\ell d_j x c_j
\]
\end{enumerate}
\end{theorem}

\begin{proof}
In the case when $H$ is finite dimensional, the statements are all
true about every $T \in \El(\BH)$. So we assume that $H$ is infinite
dimensional.

(\ref{wcpcti}) 
$\Rightarrow$
(\ref{wcpctii}) is clear by restriction of $T$.

(\ref{wcpctii}) 
$\Rightarrow$
(\ref{wcpcti}) is clear because the elementary operator on $\BH$ is
the double transpose of the operator on $\mathcal{K}(H)$.

(\ref{wcpcti}) 
$\Rightarrow$
(\ref{wcpctiii}): 
If $f_T$ is not continuous at $(0,0)$ there exist weakly null nets
$(\xi_\alpha)_\alpha$ and $(\eta_\alpha)_\alpha$ in $H_1$ so that
$f_T(\xi_\alpha, \eta_\alpha)$ does not converge to 0. Taking a subnet
we can assume that $\|\xi_\alpha\|$ and $\|\eta_\alpha\|$ are both
bounded away from 0 and then we can normalise them to be unit vectors.
Passing to
subnets again we assume $\lim_\alpha f_T(\xi_\alpha, \eta_\alpha)$ exists
and is nonzero. Passing to further subnets we can assume that
$\omega_{\xi_\alpha} \to \phi_1$ and $\omega_{\eta_\alpha} \to \phi_2$
for states $\phi_1$ and $\phi_2$ of $\BH$ vanishing on
$\mathcal{K}(H)$. Now 
\[
\tgm(Q(\mathbf{a}^*, \phi_1), Q(\mathbf{b},
\phi_2) ) \neq 0.
\]
The state $\phi = (\phi_1 + \phi_2)/2$ satisfies
\[
\tgm(Q(\mathbf{a}^*, \phi), Q(\mathbf{b}, \phi) )
\geq
\frac{1}{2}\tgm(Q(\mathbf{a}^*, \phi_1), Q(\mathbf{b}, \phi_2) )
> 0.
\]
As all states on the Calkin algebra are weak*-limits of pure states
(by \cite[11.2.4]{DixCstar} and \cite[Th.~3.3]{Gramsch}),
there must be a pure state $\psi$ of $\BH/\mathcal{K}(H)$ so that
$\tgm(Q(\mathbf{a}^*, \psi), Q(\mathbf{b}, \psi)) \neq 0$.
Thus the operator induced by $T$ on the Calkin algebra has nonzero
norm by
Theorem~\ref{CStarFormula}.

Consider an arbitrary $x \in \BH$ (of norm $\|x\| \leq 1$).
Take the net of all finite rank projections $P$ on $H$ ordered by
range inclusion. Then the net $x_P = (1-P)x(1-P)$ converges to 0 in
the weak*-topology on $\BH$ and so $T(x_P) \to 0$ in the weak*-topology (by
weak*-continuity of $T \in \El(\BH)$). By weak compactness of $T$, a
subnet of $T(x_P)$ must converge weakly to 0.

Let $\pi = \pi_\psi \colon \BH/\mathcal{K}(H) \to \mathcal{B}(H_\pi)$
be the irreducible representation
determined by the pure state $\psi$ and let
$\theta_\psi$ be the cyclic vector for the representation.
By the proofs of
Lemma~\ref{firstlemma} and Theorem~\ref{thmformulaA},
we can see that $\tgm(Q(\mathbf{a}^*, \psi),
Q(\mathbf{b}, \psi))$ is the norm of the linear functional $y \mapsto
\langle (T^\pi y) \theta_\psi, \theta_\psi \rangle$. Choose
$y \in \mathcal{B}(H_\pi)$ so that $\langle (T^\pi y) \theta_\psi,
\theta_\psi \rangle \neq 0$
and $x \in \BH$ so that $\pi(x) = y$. Then $\psi(Tx) \neq 0$.
Since $x - x_P \in
\mathcal{K}(H)$, it follows that $Tx - Tx_P \in \mathcal{K}(H)$ and so
$
|\psi(Tx_P)| = |\psi(Tx)|
> 0.
$
This contradicts any subnet $Tx_P \to 0$ weakly.

(\ref{wcpctiii}) 
$\Rightarrow$
(\ref{wcpctiv}): 
is clear.

(\ref{wcpctiv}) 
$\Rightarrow$
(\ref{wcpcti}): 
Consider a net $(x_\alpha)_\alpha$ in the unit ball of $\BH$. It has a
subnet converging weak* to some $x_0$. We call the subnet
$(x_\alpha)_\alpha$ again and aim to show $Tx_\alpha \to Tx_0$
weakly.

We denote the dual spaces of $\mathcal{K}(H)$ by $\mathcal{K}(H)'$ and
of $\BH$ by $\BH'$.
Since $\mathcal{K}(H)$ is an $M$-ideal in $\BH$,
$\BH'$ is an $\ell^1$ direct sum of those functionals vanishing on
$\mathcal{K}(H)$ (or $(\BH/\mathcal{K}(H))'$)
and a complement which is the canonical image
of $\mathcal{K}(H)'$ in its double dual $\mathcal{K}(H)''' = \BH'$.
Thus every functional $\gamma$ on $\BH$ is a
sum of (singular and normal)
functionals $\gamma_s \in (\BH/\mathcal{K}(H))'$
and $\gamma_n \in \mathcal{K}(H)'$. We know $\lim_\alpha
\gamma_n(Tx_\alpha) = \gamma_n(Tx_0)$ by weak*-continuity of $T$ and so we
concentrate on establishing $\lim_\alpha 
\gamma_s(Tx_\alpha) = \gamma_s(Tx_0)$.

As $\gamma_s$ can be expressed as a linear combination of 4 states, it
is enough to deal with the case when $\gamma_s = \phi$ is a state of the
Calkin algebra.
But in that case, we know we can express $\phi$ as a weak*-limit
$\phi(y) = \lim_\beta \omega_{\xi_\beta}(y)$ ($y \in \BH$) for
a weakly null net $(\xi_\beta)_\beta$ of unit vectors in $H$.
The norm of the functional $x \mapsto \phi(Tx)$ is
\begin{eqnarray*}
\tgm(Q(\mathbf{a}^*, \phi), Q(\mathbf{b}, \phi) )
&=&
\lim_\beta \tgm(Q(\mathbf{a}^*, \omega_{\xi_\beta}), Q(\mathbf{b}, \omega_{\xi_\beta}) )\\
&=&
\lim_\beta f_T(\xi_\beta, \xi_\beta) = 0
\end{eqnarray*}
by (\ref{wcpctiv}). Hence $\phi(Tx_\alpha) = \phi(Tx_0) = 0$
for all $\alpha$.

(\ref{wcpctiii}) 
$\Rightarrow$
(\ref{wcpctv}): 
Consider a maximal subset of $\{a_1, a_2, \ldots, a_\ell \}$ which is
linearly independent modulo $\mathcal{K}(H)$. By renumbering, we can
assume this maximal set is $a_1, a_2, \ldots, a_m$ for $0 \leq m \leq
\ell$.  (If $m=0$, all the $a_j$ are compact.) We can then express
$a_j$ for $m+1 \leq j \leq \ell$ as a compact $c_j$ plus a linear
combination of $a_1, a_2, \ldots, a_m$. This allows us to write $T$ in
the form
\[
Tx = \sum_{j=1}^m a_j x b_j' + \sum_{j=m+1}^\ell c_j x b_j = T_1x +
T_2x.
\]
It is easy to see that $f_{T_2}$ is continuous at points $\{0\} \times
H_1$. As $f_T(\xi, \eta)$ is the norm of the functional $x \mapsto \langle
(Tx) \eta, \xi \rangle$, it follows that $f_{T_1}(\eta, \xi) \leq f_T(\eta,
\xi) +
f_{T_2}(\eta, \xi)$ and is continuous at $(0,0) \in H_1 \times
H_1$.

Thus it is sufficient to consider the case where $m= \ell$ and the
$a_j$ are linearly independent modulo $\mathcal{K}(H)$ and to 
show that each $b_j$ is compact in this case.
By Lemma~\ref{EssLinIndepLemma} there is a weakly null net
$(\xi_\alpha)_\alpha$ in $H_1$ so that
$\lim_\alpha Q(\mathbf{a}^*, \xi_\alpha)$ is positive definite. Thus
for $\alpha$ large there is $\epsilon >0$ so that $Q(\mathbf{a}^*,
\xi_\alpha) > \varepsilon I_\ell$. For any net $(\eta_\beta)_\beta$ in
$H_1$ which is weakly null, we can order pairs $(\alpha, \beta)$ via
$(\alpha_1, \beta_1) \leq (\alpha_2, \beta_2) \iff \alpha_1 \leq
\alpha_2$ and $\beta_1 \leq \beta_2$ and thereby create a net
$((\xi_\alpha, \eta_\beta))_{(\alpha, \beta)}$ which converges to
$(0,0)$ in $H_1 \times H_1$. However
\[
f_T(\xi_\alpha, \eta_\beta) = \tgm( Q(\mathbf{a}^*, \xi_\alpha),
Q(\mathbf{b}, \eta_\beta) ) \geq \sqrt{\varepsilon} \sqrt{
\trace Q(\mathbf{b},
\eta_\beta)}
\]
and it follows that $\lim_\beta \|b_j\eta_\beta\| = 0$ for each $j$.
Thus each $b_j$ is compact in this case.

(\ref{wcpctv}) 
$\Rightarrow$
(\ref{wcpctiii}): 
is easy to verify by writing
\[
Tx = \sum_{j=1}^m c_j x d_j + \sum_{j=m+1}^\ell d_j x c_j = T_1 x +
T_2x,
\]
and using $f_T \leq f_{T_1} + f_{T_2}$.
It is easy to show show that $f_{T_1}$ is continuous at points of
$\{0\} \times H_1$ and $f_{T_2}$  is continuous at points of $H_1
\times \{0\} $.
\end{proof}

A natural question which remains unresolved is whether an
analogue of Theorem~\ref{BHweakcompactness} holds for weakly compact
elementary operators on (general) $C^*$-algebras. There are several
related results established in \cite[\S5.3]{AraMatBk}.

\vfill

\noindent
\textbf{Keywords:} Tracial geometric mean; matrix numerical range;
$C^*$-algebra\\
2000 Mathematics Subject Classification: 47B47, 46L07
\vfill

\noindent
School of Mathematics\\
Trinity College\\
Dublin 2\\
Ireland\\[2mm]
Email: \texttt{richardt@maths.tcd.ie}

\vfill

\begin{thebibliography}{00}
\bibitem{JCAlexander} J. C. Alexander, \emph{Compact Banach
algebras}, Proc. London Math. Soc. (3) \textbf{18} (1968), 1--18.
\bibitem{AndoEtAl} T. Ando, C.-K. Li\ and\ R. Mathias,
\emph{Geometric means}, Linear Alg. Appl. \textbf{383} (2004) 305--334.
\bibitem{AraMatBk} P. Ara\ and\ M. Mathieu,  Local multipliers of
$C$*-algebras, Springer, London (2003).
\bibitem{Arch78}
R. J. Archbold, \emph{On the norm of an inner derivation of a
$C\sp*$-algebra}, Math. Proc. Cambridge Philos. Soc. \textbf{84} (1978),
273--291.
\bibitem{AMS} R. J. Archbold, M. Mathieu\ and\ D. W. B. Somerset,
\emph{Elementary operators on antiliminal
$C$*-algebras}, Math. Ann. \textbf{313} (1999), 609--616.
\bibitem{ArchSomerset04}
R. J. Archbold\ and\ D. W. B. Somerset, \emph{Inner derivations and
primal ideals of $C\sp *$-algebras. II}, Proc. London Math. Soc. (3)
\textbf{88} (2004), 225--250
\bibitem{AST} R. J.~Archbold, D. W.~B.~Somerset and R. M. Timoney,
\emph{On the central Haagerup
tensor product and completely bounded
mappings of a $C^*$-algebra}, J. Funct. Anal. (to appear).
\bibitem{ArchBatty} C.~J.~K.~Batty and R.~J.~Archbold, On
factorial states of operator algebras II, \emph{J. Operator Theory}
\textbf{13} (1985), 131--142.
\bibitem{Bhatia} R. Bhatia, Matrix Analysis, Graduate Texts in
Mathematics 169, Springer-Verlag, New York (1997).
\bibitem{BBR}
A. Blanco, M. Boumazgour\ and\ T. J. Ransford, \emph{On
the norm of elementary operators}, J. London Math. Soc.
(2) {\bf 70} (2004), 479--498

\bibitem{DixCstar} J. Dixmier, Let $C^*$-alg\`ebres et leurs
representations, Gauthier-Villars, Paris 1964.
\bibitem{EffrosRuan} E. G. Effros\ and\ Z.-J. Ruan, Operator
Spaces, London Mathematical Society Monographs 23, Oxford Science
Publications (2000).
\bibitem{FongSourour} C. K. Fong and A. R. Sourour,
\emph{On the operator identity $\sum \,A\sb{k}XB\sb{k}\equiv 0$},
Canad. J. Math. \textbf{31} (1979), 845--857.
\bibitem{Gramsch}
B. Gramsch, 
\emph{Eine Idealstruktur Banachscher Operatoralgebren},
J. Reine Angew. Math. \textbf{225} (1967), 97--115
\bibitem{KaplanskyNormed} I. Kaplansky, \emph{Normed algebras}, Duke
Math. J. \textbf{16} (1949) 399--418.
\bibitem{Magajna93}
B. Magajna,
 \emph{A transitivity theorem for algebras of elementary operators},
 Proc. Amer. Math. Soc. \textbf{118} (1993), 119--127.
\bibitem{Magajna} B. Magajna,
\emph{The Haagerup norm on the tensor
product of operator modules}, J. Funct. Anal. \textbf{129} (1995),
no. 2, 325--348.
\bibitem{MagTurn} B. Magajna\ and\ A. Turn\v sek,
\emph{On the
norm
of symmetrised two-sided multiplications}, Bull. Australian Math. Soc.
\textbf{67} (2003) 27--38.
\bibitem{MMNormProb} M. Mathieu, \emph{The norm problem for
elementary operators}, in Recent Progress in Functional Analysis (K.
D. Bierstedt et al, eds.), Elsevier, Amsterdam 2001, pp. 363--368.
\bibitem{Paulsen}
V. Paulsen, \emph{Completely bounded maps and operator algebras},
Cambridge Univ. Press, Cambridge (2002).
\bibitem{Somerset93}
D. W. B. Somerset, \emph{The inner derivations and the primitive
ideal space of a $C\sp *$-algebra}, 
J. Operator Theory \textbf{29} (1993), 307--321.
\bibitem{Somerset94}
D. W. B. Somerset, \emph{Inner derivations and primal ideals of $C\sp
*$-algebras}, J. London Math. Soc. (2) {\bf 50} (1994),
568--580
\bibitem{Somerset97}
D. W. B. Somerset, \emph{The proximinality of the centre of a $C\sp
*$-algebra}, J. Approx. Theory {\bf 89} (1997), 114--117
\bibitem{Somerset98}
D.~W.~B. Somerset, \emph{The central Haagerup tensor
product of a $C$*-algebra}, J. Operator Theory \textbf{39} (1998)
113--121.
\bibitem{Stampfli} J. G. Stampfli, \emph{The norm of a
derivation}, Pacific J. Math. \textbf{33} (1970) 737--747.
\bibitem{TimoneyNE} R. M. Timoney, \emph{Computing the norms of
elementary operators}, Illinois J. Math \textbf{47} (2003) 1207--1226.
\bibitem{Timoneyaxb} R. M. Timoney,
\emph{Norms and CB norms of Jordan elementary operators}, 
Bull. Sci. Math. {\bf 127} (2003), 597--609.
\bibitem{Vala} K. Vala, \emph{On compact sets of compact operators},
Ann. Acad. Sci. Fenn. Ser A I \textbf{351} (1964).
\end{thebibliography}
\end{document}